\documentclass[12pt]{amsart}
\usepackage{amssymb}
\usepackage{amsmath, amsthm, amscd, amsfonts, amssymb, graphicx,mathrsfs, color}
\usepackage[colorlinks=true, linkcolor=blue,urlcolor=blue]{hyperref}
\input xy \xyoption{all}
\usepackage[all]{xy}
\usepackage{layout}
\usepackage{tikz}
\usepackage{bm}
\usepackage{geometry}
\usepackage{euscript}
\usepackage{enumerate}

\newtheorem{thm}{Theorem}[section]
\newtheorem{cor}[thm]{Corollary}
\newtheorem{lem}[thm]{Lemma}
\newtheorem{prop}[thm]{Proposition}
\theoremstyle{definition}
\newtheorem{defn}[thm]{Definition}
\newtheorem{ex}[thm]{Example}
\theoremstyle{remark}
\newtheorem{rem}[thm]{Remark}

\numberwithin{equation}{section}

\begin{document}

\author[A.N. Abyzov]{Adel N. Abyzov}
\address{Department of Algebra and Mathematical Logic, Kasan Federal University, Volga Region, Russian Federation}
\email{adel.abyzov@kpfu.ru}

\author[Ruhollah Barati]{Ruhollah Barati}
\address{Department of Pure Mathematics, Faculty of Mathematical Sciences,
Tarbiat Modares University, Tehran, Iran}
\email{baratiroohollah@gmail.com; ruhollahbarati@modares.ac.ir}

\author[P.V. Danchev]{Peter V. Danchev$^*$}
\address{Institute of Mathematics and Informatics, Bulgarian Academy of Sciences, Sofia, Bulgaria}
\email{danchev@math.bas.bg; pvdanchev@yahoo.com}

\title{Rings Close to Periodic with Applications to Matrix, Endomorphism and Group Rings}
\thanks{*Corresponding author}

\begin{abstract} We examine those matrix rings whose entries lie in periodic rings equipped with some additional properties. Specifically, we prove that the famous Diesl's question whether or not {\it $R$ being nil-clean implies that $\mathbb{M}_n(R)$ is nil-clean for all $n\geq 1$} is paralleling to the corresponding implication for (Abelian, local) periodic rings. Besides, we study when the endomorphism ring $\mathrm{E}(G)$ of an Abelian group $G$ is periodic. Concretely, we establish that $\mathrm{E}(G)$ is periodic exactly when $G$ is finite as well as we find a complete necessary and sufficient condition when the endomorphism ring over an Abelian group is strongly $m$-nil clean for some natural number $m$ thus refining an "old" result concerning strongly nil-clean endomorphism rings. Responding to a question when a group ring is periodic, we show that if $R$ is a right (resp., left) perfect periodic ring and $G$ is a locally finite group, then the group ring $RG$ is periodic, too. We finally find some criteria under certain conditions when the tensor product of two periodic algebras over a commutative ring is again periodic. In addition, some other sorts of rings very close to periodic rings, namely the so-called weakly periodic rings, are also investigated.
\end{abstract}

\maketitle

\section*{Introduction and Motivation}\label{IntroductionAndMotivation}

Throughout the text of the current article, all rings are associative with an identity. Almost all notions and notations are standard being in agreement with the well-known book \cite{L}, and the more specific terminology will be stated explicitly in what follows.

We denote the Jacobson radical, the lower nil-radical, the upper nil-radical, the Levitzki radical, the set nilpotent elements, and the center of a ring $R$,  by the symbols $J(R)$, $\mathrm{Nil}_*(R)$, $\mathrm{Nil}^*(R)$, $L\text{-}rad(R)$, $\mathrm{Nil}(R)$, and $\mathrm{C}(R)$, respectively. Moreover, note that the following inclusions hold:

\[\mathrm{Nil}_*(R)\subseteq L\text{-}rad(R)\subseteq \mathrm{Nil}^*(R)\subseteq \mathrm{Nil}(R)\cap J(R).\]

\medskip

Recall that a ring $R$ is said to be {\it $2$-primal} if $\mathrm{Nil}_*(R)=\mathrm{Nil}(R)$. For example, every commutative or reduced ring is $2$-primal. Imitating \cite{wei-xing}, the ring $R$ is called {\it weakly $2$-primal} if the equality $\mathrm{Nil}(R)=L\text{-}rad(R)$ holds. It is clear that every $2$-primal ring is weakly $2$-primal. However, the converse implication fails as there exist weakly $2$-primal rings that are not $2$-primal (see, for instance, \cite[Example 2.2]{Marks}). Also, a ring $R$ is called $\mathrm{NI}$ if the set $\mathrm{Nil}(R)$ forms an ideal.

It was shown in \cite[Lemma 3.12]{BMA} that if $R$ is an $m$-potent ring for some arbitrary but fixed integer $m>1$, that is $x^m=x$ for all $x\in R$, then $\mathbb{M}_n(R)$ is a periodic ring for all naturals $n\geq 1$. This assertion was considerably generalized in \cite[Corollary 2.22]{BCL} for any non-fixed $m\geq 2$ depending on the element $x$. Some other important results in this aspect are \cite[Proposition 2.18, Theorem 2.19, Theorem 2.21]{BCL}. On the same vein, in \cite{D1} the "periodicity" was extended to the concept of so-called {\it $\pi$-UU rings}; recall that a ring $R$ is said to be $\pi$-UU if, for each invertible element $u\in U(R)$, the unit group of $R$, there exists an integer $i>0$ which depends on $u$ such that $u^i\in 1+Nil(R)$ (for the case of fixed $i=1$, we refer to \cite{DL} calling them just {\it UU rings}). It is clear that all periodic rings are immediately $\pi$-UU rings, that implication is obviously non-reversible; however, it was proved in \cite[Theorem 3.4 (7)]{CD} that strongly $\pi$-regular $\pi$-UU rings are themselves periodic -- we refer the interested reader to the good sources \cite{G} and \cite{T} for more information and properties concerning the strong $\pi$-regularity of rings.

On the other hand, in \cite{D13} and \cite{BDZ} (compare also with \cite{DL} and \cite{KWZ}) were considered some versions of the classical cleanness of rings, originated by W.K. Nicholson in \cite{N}, such as {\it nil-cleanness}, {\it strong nil-cleanness} and {\it weak nil-cleanness}, respectively. It is worthwhile noticing that strongly nil-clean rings are always strongly $\pi$-regular, where the latter concept is defined as in \cite{L}. Some non-trivial examples of such rings can be found in \cite{HKL} and the bibliography cited therewith. Some related papers of this topic are \cite{BMA} and \cite{Di} as well, which deal with rings whose elements are sums of commutating potent and nilpotent elements.

Likewise, it was proved in \cite[Corollary 3.10]{BMA} that the matrix ring $\mathbb{M}_n(R)$ is always nil-clean for all $n\geq 1$, provided that the former ring $R$ is weakly $2$-primal strongly nil-clean -- actually, the authors stated this claim without the assumption of being "weakly $2$-primal" by using \cite[Proposition 2.7]{NI} which proposition seems to be false. (For the partial case of a $2$-primal ring, we refer the interested readers for more information to \cite{D2} and to \cite{BCDM} as well). Moreover, a similar topic to that in \cite{BMA} was considered in \cite{PDK}, but unfortunately the proof of \cite[Lemma 3.3]{PDK} is erroneous when the integer $m\geq 2$ is {\it not} a prime -- in fact, the ring $\mathbb{Z}_4$ manifestly demonstrates our claim as simple computations show. However, certain appropriate results from \cite{khurana2021lifting} could be applied to avoid these shortcomings. Indeed, applying \cite[Theorem 8]{khurana2021lifting}, {\it letting $R$ be a ring in which $m-1$ is a unit, then $m$-potents lift modulo every nil ideal $N$ of $R$}.

Unifying these two directions, it was asked in \cite[p.9, Question]{CD} of whether or not nil-clean $\pi$-UU rings are always periodic. The purpose of this paper is to present some further results and relationships of this branch by examining periodic rings and their close versions in a rather more attractive manner as follows. Concretely, our further work is organized in the sequel thus:

In the first section, we find a satisfactory criterion when an arbitrary group ring is periodic. We succeeded here to establish some substantial affirmations (see Theorems~\ref{theorem3.4}, \ref{theorem3.7}). Moreover, we also explore when a group ring is periodic by answering in part a question posed by Danchev in Mat. Stud. (2020). Precisely, we show that if $R$ is a right (resp., left) perfect periodic ring and $G$ is a locally finite group, then the group ring $RG$ is periodic, too.

The objective of the second one is to study the behavior of the matrix rings over a periodic ring. Our achievements here are Theorems~\ref{theorem2.3}, \ref{theorem2.5}, \ref{theorem2.7}. Specifically, we prove that the full $n\times n$ matrix ring $\mathbb{M}_n(R)$ is periodic for all natural numbers $n$, whenever the base ring $R$ is weakly $2$-primal periodic. This somewhat substantially refines recent results in the subject established by Cui-Danchev in J. Algebra \& Appl. (2020) and by Bouzidi-Cherchem-Leroy in Commun. Algebra (2021).

The leitmotif of the third one is to find necessary and sufficient conditions when an arbitrary endomorphism ring over an Abelian group is either periodic or strongly $m$-nil clean for some integer $m\geq 1$ (notice that the latter class of rings is always contained in the class of periodic rings). We here achieved to prove that the endomorphism ring is periodic exactly when the former Abelian group is finite as well as that such a ring is strongly $m$-nil clean uniquely when the Abelian group is finite and its $p$-primary components have some special properties (cf. Theorem~\ref{theorem3.11}, \ref{theorem3.12} and Corollary \ref{corollary3.13}). Thereby, we continue and strengthen the well-known earlier results due to Fuchs-Rangaswamy in Math. Z. (1968) about $\pi$-regular endomorphism rings as well as we refine a recent result due to Breaz-C\v{a}lug\v{a}reanu-Danchev-Micu in Lin. Algebra \& Appl. (2013) for strongly nil-clean endomorphism rings.

In the fourth one, we examine the tensor products of algebras in the context of the periodicity as our basic results are, respectively, structured in Theorems~\ref{theorem5.2} and \ref{theorem5.4}. Thus, our results supply those from \cite{S}.

In the fifth section, we consider the so-called {\it weakly periodic rings} and present some basic properties of them. Our motivating tool to do that is the fact that these rings are situated between the well-known classes of nil-clean rings in \cite{D13} and semi-clean rings in \cite{Ye}, as well as that the class of weakly periodic rings properly contains the well-known class of periodic rings as it will be shown below. Our results here are, respectively, presented in Theorems~\ref{theorem5.5} and \ref{theorem3.10}, which can be viewed as a natural continuation of those obtained by Cui-Danchev in J. Algebra \& Appl. (2020).

We close our work in the sixth section by posing six well-arranged questions of some interest and importance which, hopefully, will motivate a possible further intensive research on the explored object.

In the spirit of all given sections, we will try to demonstrate some logical relationship between some of them in light of the received results.

\section{Group Rings} \label{SectionGroupRing}

\medskip

It was established in \cite[Theorem 1]{DanchevGroupRings} that the commutative group ring $RG$ is periodic if and only if $R$ is a periodic ring and $G$ is a torsion group, provided $R$ is a commutative local ring with $p\cdot 1\in \mathrm{Nil}(R)$ for some prime $p$ and $G$ is an Abelian group.

Before attempting to establish our chief results in this section, which will improve the cited one, we need some more conventions in what follows.

First, we recall once again some well-known definitions like these: A ring $R$ (respectively, a group $G$) is said to be {\it locally finite} if any finitely generated subring of $R$ (respectively, any finitely generated subgroup of $G$) is a finite ring (respectively, a finite group). Also, a ring $R$ is said to be of {\it bounded index} of nilpotence if there is a positive integer $n$ such that $a^n=0$ for any nilpotent element $a$ in $R$. Finally, a ring $R$ is said to be of {\it locally bounded index} if every finitely generated subring of $R$ is of bounded index.

According to \cite[Lemma 2.1]{BCL}, if $R$ is a periodic ring, then its characteristic is always positive. Therefore, we will restate \cite[Corollary 2]{Hirano} as follows.

\begin{lem}[{\cite[Corollary 2]{Hirano}}]\label{lemma3.1}
A ring $R$ is locally finite if and only if $R$ is a periodic ring of locally bounded index.
\end{lem}

We are now in a position to state and prove our basic assertions. Before doing that, we just note that (*) {\it if the group ring $RG$ is periodic, then the former ring $R$ is periodic too as its epimorphic image}. But, moreover, being periodic, the group ring $RG$ is known to be strongly $\pi$-regular, so we extract that $G$ is a torsion group by employing \cite[Proposition 3.4]{HV}. In the following, we will prove the statement (*) with a different approach. Also, we shall prove the converse of this statement with some additional assumptions on both $R$ and $G$.

\begin{thm}\label{theorem3.4}
If $R$ is a weakly $2$-primal periodic ring and $G$ is a locally finite group, then the group ring $RG$ is periodic.
\end{thm}

\begin{proof}
Take an arbitrary element $\sum_{i=1}^n a_ig_i \in \mathrm{Nil}(R)G$. As $\mathrm{Nil}(R)$ is locally nilpotent, for the finite set of elements $\{a_1,a_2, \cdots ,a_n\}$ there exists a natural number $N$ such that any product of $N$ elements from $\{a_1,a_2, \cdots ,a_n\}$ is zero. But ${(\sum_{i=1}^n a_ig_i)}^N$ is a sum of the items $a_{i_1}a_{i_2} \cdots a_{i_N}g_{i_1}g_{i_2} \cdots g_{i_N}$. So, it follows that ${(\sum_{i=1}^n a_ig_i)}^N=0$. Therefore, $\mathrm{Nil}(R)G$ is a nil-ideal of $RG$. Since the quotient $R/\mathrm{Nil}(R)$ is a potent ring, the application of Lemma~\ref{lemma3.1} yields that it is a locally finite ring. However, by \cite[Proposition 2.12]{HKL}, the group ring $(R/\mathrm{Nil}(R))G$ is locally finite, too. Now, by combining the isomorphism $\frac{RG}{\mathrm{Nil(R)}G} \cong (\frac{R}{\mathrm{Nil}(R)})G$ with Lemma~\ref{lemma3.1} and \cite[Corollary 3.6]{CD}, we can conclude the desired result.
\end{proof}

As two immediate consequences of the last assertion, we obtain the following.

\begin{cor}\label{corollary3.5}
If $R$ is a $2$-primal periodic ring and $G$ is a locally finite group, then the group ring $RG$ is periodic.
\end{cor}

\begin{cor}\label{corollary3.6}
Let $R$ be a commutative ring and $G$ an Abelian group. Then, the group ring $RG$ is periodic if and only if $R$ is periodic and $G$ is torsion.
\end{cor}

A possible non-trivial extension of the last result, which we are unable to prove or disprove at present, is the following: If $R$ is a periodic ring such that $J(R)$ is locally nilpotent and $G$ is a locally finite group, then the group ring $RG$ is periodic. In fact, if it is true that $R$ contains no finite subset of non-zero orthogonal pairwise-isomorphic idempotents, then the quotient-ring $R/J(R)$ will have a bounded index of nilpotency by \cite[Theorem 2.8]{W}. Now, the group ring $(R/J(R))G$ will be locally finite by Lemma~\ref{lemma3.1} and
\cite[Proposition 2.12]{HKL}. Hence, $(R/J(R))G$ has to be periodic. Moreover, it is well known that $J(R)G$ is nil.
Thus, the isomorphism $RG/J(R)G \cong (R/J(R))G$ will imply the desired result.

\medskip

In \cite{P} and \cite{DM}, respectively, were found criteria when an arbitrary ring is potent in terms of sub-direct products of finite fields. In what follows, we shall completely characterize when the group ring $RG$ is potent. By mimicking \cite{passman}, let $\pi(R)=\{p_1,p_2, \cdots,p_k \}$ be the set of all prime divisors of the positive characteristic of $R$. Thus, we have the following theorem.

\begin{thm}
The group ring $RG$ is potent if and only if the ring $R$ is potent and $G$ is a torsion Abelian group such that the order of any finite subgroup of $G$ is not divisible by any $p \in \pi(R)$.
\end{thm}

\begin{proof}
($\Rightarrow$). Suppose that the group ring $RG$ is potent. As $R$ is a subring of $RG$, it is too a potent ring and its characteristic is positive. Moreover, it follows from Corollary~\ref{corollary3.6} that $G$ is a torsion Abelian group. If now $p \in \pi(R)$, then by \cite[Theorem III]{passman} the group $G$ has no a finite subgroup whose order is divisible by $p$.

($\Leftarrow$). Consulting with Corollary~\ref{corollary3.6}, the group ring $RG$ is periodic and, besides, it is simultaneously reduced according to \cite[Theorem I]{passman}. Therefore, with the aid of \cite[Theorem 3.4 (1),(2)]{CD}, we can derive that the group ring $RG$ is potent, as claimed.
\end{proof}

Recall that a ring is said to be {\it right perfect} if the ideal $J(R)$ vanishes on the right (that is, for any sequence $a_1, a_2, \cdots ,$ of elements of $J(R)$, there is an integer, say $k$, such that their product is zero, i.e., $a_1\cdots a_k =0$) and the quotient-ring $R/J(R)$ is semi-simple.

\begin{thm}\label{theorem3.7}
Let $R$ be a right (resp., left) perfect periodic ring and let $G$ be a locally finite group. Then, $RG$ is a periodic ring.
\end{thm}

\begin{proof}
If $R$ is a right (resp., left) perfect ring, then the fact that $J(R)$ is locally nilpotent is guaranteed by \cite[Proposition 23.15]{L} (see also \cite[Exercise 23.1]{L}). Now, similarly to the proof of Theorem \ref{theorem3.4} alluded to above, we derive that $J(R)G$ is nil. Since $R/J(R)$ is semi-simple, by the help of the old-standing famous Artin-Wedderburn theorem, we may write that $R/J(R) \cong \prod_{i=1}^k \mathbb{M}_{n_i}(D_i)$ for some division rings $D_i$. Being a periodic ring, the factor-ring $R/J(R)$ implies that each $D_i$ is also periodic and, consequently, it is a locally finite field owing to Lemma \ref{lemma3.1}. Therefore, $\mathbb{M}_{n_i}(D_i)$ is a locally finite ring applying \cite[Corollary 2.3]{HKL}. Now, it follows from \cite[Proposition 2.12]{HKL} that each of the group rings $[\mathbb{M}_{n_i}(D_i)]G$ over the ring $\mathbb{M}_{n_i}(D_i)$ is locally finite and, consequently, it is periodic employing Lemma \ref{lemma3.1}. Furthermore, the isomorphism $(\frac{R}{J(R)})G \cong \prod_{i=1}^k \mathbb{M}_{n_i}(D_i)G$ and \cite[Lemma 2.12]{BCL} yield that $\frac{RG}{J(R)G} \cong (\frac{R}{J(R)})G$ is a periodic ring, as well. Finally, an application of \cite[Corollary 3.6]{CD} gives the desired result.
\end{proof}

\begin{rem}
Note that Theorem~\ref{theorem3.4} and Theorem~\ref{theorem3.7} give only partial answers to \cite[Problem 1]{DanchevGroupRings}.
\end{rem}

\section{Matrix Rings}\label{SectionMatrixRings}

\medskip

Referring to \cite[Corollary 2.25]{BCL} (see \cite{CD} too), one can see that if $R$ is a ring whose elements satisfy the equation $x^m=x^n$, where $m> n$ are two fixed natural numbers, then the full matrix ring $\mathbb{M}_k(R)$ remains periodic for any integer $k\geq 1$. We shall now characterize these rings $R$ as follows (see \cite[Proposition 4]{D3} as well).

\begin{prop} \label{proposition2.2}
For a ring $R$ and given two natural numbers $m > n$, the followings two conditions are equivalent:
\begin{enumerate}
\item[\rm{(1)}] $x^m=x^n$ for all $x \in R$.

\item[\rm{(2)}] In a ring $R$, every element $x$ can be represented as $x=a+b$, where $b$ is a $(m-n+1)$-potent,  $a^n=0$, $ab=ba=0$ and the characteristic of $R$ is finite.
\end{enumerate}
\end{prop}

\begin{proof}
(1) $\Rightarrow$ (2). Given any non-negative integers $m > n$, we have the following isomorphism
$$\mathbb{Z}[x]/(x^m-x^n) \cong \mathbb{Z}[x]/(x^n) \times \mathbb{Z}[x]/(x^{m-n}-1)$$
utilizing the Chinese Remainder Theorem. Thus, the element $x$ automatically decomposes as a sum of a nilpotent (of index no greater than $n$) and an $(m-n + 1)$-potent that commute each other. In fact, their left/right products are zero. Moreover, when the equality $x^m = x^n$ holds for every $x \in R$ (with $m$ and $n$ possibly depending on $x$), then this situation forces finite characteristic $d$ of $R$ just by taking $x = 2$.

(2) $\Rightarrow$ (1). Given an element $x = a+b$, where $ab = ba = 0$, $a^n=0$ and $b^{m-n+1} = b$, we quickly obtain by a direct calculation that $x^m = x^n$.
\end{proof}

\begin{rem}
\begin{enumerate}
\item[\rm{(1)}]
In Proposition \ref{proposition2.2}, when $m = n + 2$, we get
$$\mathbb{Z}_d[x]/(x^{m-n}-1) = \mathbb{Z}_d[x]/(x^2-1) \cong \mathbb{Z}_d[x]/(x-1) \times \mathbb{Z}_d[x]/(x+1) \cong \mathbb{Z}_d \times \mathbb{Z}_d.$$

\noindent So, similar arguments should apply here, too. We, thus, infer that $d$ is a divisor of $2^{n+2}-2^n = 2^n · 3$.
\item[\rm{(2)}]
In Proposition \ref{proposition2.2}, if $x^m=x^n$ for all $x \in R$ and $m$ and $n$ have opposite partiy then $R$ is a potent ring by \cite[Corollary 2.2]{Anderson}.
\end{enumerate} 
\end{rem}

Call a ring {\it right (resp., left) quasi-duo}, provided every maximal right (resp., left) ideal is two-sided (see, e.g., \cite{yu1995quasi}).
Weakly periodic rings are introduced in section \ref{SectionWeaklyPeriodicRings}. In the next proposition, if $R$ is a weakly periodic then the statement holds.

\begin{prop}\label{proposition2.1}
For a ring $R$, the following three statements are equivalent.
\begin{enumerate}
\item[\rm{(1)}]
$J(R)$ is nil and $R/J(R)$ is a potent ring.
\item[\rm{(2)}]
$R$ is a right (left) quasi-duo (weakly )periodic ring.
\item[\rm{(3)}]
$R$ is a (weakly) periodic NI-ring.
\end{enumerate}
\end{prop}

\begin{proof}
It is straightforward. We just need to show the implication (2) $\Rightarrow$ (1). It is clear that $J(R)$ is nil. Since $R$ is right (left) quasi-duo, it follows that the factor-ring $R/J(R)$ is reduced in view of \cite[Corollary 2.4]{yu1995quasi}. Therefore, the quotient $R/J(R)$ is potent, as needed.
\end{proof}

The next result could somewhat be viewed as an generalization of \cite[Theorem 2.6]{DGL}.

\begin{thm}\label{theorem2.3}
Suppose that $R$ is a weakly $2$-primal periodic ring. Then, the full matrix ring $\mathbb{M}_n(R)$ is periodic for all $n \geq 1$.
\end{thm}

\begin{proof}
As $J(R)$ is nil, $J(R)=\mathrm{Nil}(R)=L\text{-}rad(R)$. Now, a simple combination between the isomorphism $\mathbb{M}_n(R)/\mathbb{M}_n(J(R)) \cong \mathbb{M}_n(R/J(R))$, Proposition \ref{proposition2.1}, \cite[Corollary 2.22]{BCL}, \cite[Proposition 1.1]{BM} and \cite[Corollary 3.6] {CD} imply the desired result, as expected.
\end{proof}

As a consequence, we directly yield the following:

\begin{cor}\label{corollary2.4}
If $R$ is a $2$-primal periodic ring, then the matrix ring $\mathbb{M}_n(R)$ is periodic for all $n\geq 1$.
\end{cor}

We now intend to consider the so-called Armendariz rings and their generalizations as follows:

\begin{defn}[{\cite[Definition 2.1]{Nasr-Isfahani}}]
Let $R$ be a ring with an endomorphism $\alpha$. We say that $R$ is a \textit{skew-Armendariz} ring, if for any two polynomials $f(x)=a_0+a_1x+a_2x^2+ \cdots +a_nx^n$ and $g(x)=b_0+b_1x+b_2x^2+ \cdots +b_mx^m$ in $R[x;\alpha]$, the equality $f(x)g(x)=0$ holds if and only if $a_ib_j=0$ for each $0 \leq i \leq n$, $0 \leq j \leq m$. In particular, if $\alpha=id_R$, then the ring $R$ is called \textit{Armendariz}.
\end{defn}

We, thus, may obtain the following result.

\begin{cor}\label{corollary2.4-1}
If $R$ is a \textit{skew-Armendariz} periodic ring, then the full matrix ring $\mathbb{M}_n(R)$ is periodic.
\end{cor}

\begin{proof}
It is well known that every \textit{skew Armendariz} ring is Abelian. But, on the other hand, each Abelian periodic ring is always an NI-ring. Now, \cite[Theorem 2.8]{Nasr-Isfahani} applies to get that the ring $R$ is $2$-primal. Therefore, Corollary \ref{corollary2.4} implies that the full matrix ring $\mathbb{M}_n(R)$ is periodic, as claimed.
\end{proof}

The next two observations are also interesting enough to be documented.

\begin{ex}
Let $p$ be a prime and $m$ an natural number. Then, for every $s$ with $[\frac{m}{2}] \leq s \leq m$, the full matrix ring $\mathbb{M}_n(\mathbb{Z}_{p^m} \propto p^s\mathbb{Z}_{p^m})$ is periodic.
\end{ex}

\begin{proof}
By \cite[Theorem 2.13]{BCL}, $\mathbb{Z}_{p^m} \propto p^s\mathbb{Z}_{p^m}$ is a periodic ring. However, it follows from \cite[Corollary 2.8]{LeeZhou} that $\mathbb{Z}_{p^m} \propto p^s\mathbb{Z}_{p^m}$ is Armendariz. Now, Corollary \ref{corollary2.4-1} gives the desired result.
\end{proof}

\begin{ex}
In \cite[Example 1.7]{BM}, if $R$ is a potent ring and $S\in \{A_n(R),A_n(R)+RE_{1,k},T(R,n),S_n^k(R)\}$, then the ring $\mathbb{M}_t(S)$ is periodic for any $t \geq 1$.
\end{ex}

We now continue with the following result of some interest and importance for our further study.

\begin{thm}\label{theorem2.5}
If $R$ is a right (resp., left) perfect periodic ring, then the full matrix ring $\mathbb{M}_n(R)$ is periodic.
\end{thm}

\begin{proof}
If $R$ is a right (resp., left) perfect ring, then $J(R)$ is locally nilpotent by \cite[Proposition 23.15]{L} (see also \cite[Ex. 23.1]{L}). But exploiting \cite[Proposition 1.1]{BM}, we know that $J(\mathbb{M}_n(R))$ is nil. The rest of our argument is similar to the one in \cite[Proposition 2.18]{BCL}.
\end{proof}

The next observations give a connection between the results from the preceding section and the current one.

\begin{rem}
If $R$ is a right (resp., left) perfect periodic ring of characteristic $2$, then by combination of Theorem \ref{theorem3.7} and the isomorphism $R\mathbb{S}_3 \cong RC_2 \oplus \mathbb{M}_2(R)$, where $\mathbb{S}_3$ is the symmetric $3$-group and $C_2$ is the cyclic group of order $2$  -- see, e.g., \cite[Lemma 4.4]{morphicgrouprings}) -- will imply that the ring $\mathbb{M}_2(R)$ is periodic.
\end{rem}

As an immediate consequence, we yield:

\begin{cor}[{\cite[Proposition 2.18]{BCL}}]\label{corollary2.6}
If $R$ is an Artinian periodic ring, then the ring $\mathbb{M}_n(R)$ is periodic for any $n \geq 1$.
\end{cor}

\medskip

We also pose here the following conjecture which is relevant to the above Theorem~\ref{theorem2.5} and whose positive solution will allow us to resolve in the affirmative the difficult problem due to Diesl concerning nil-clean matrix rings -- compare with Remark~\ref{Dis} below and, for more details, see \cite{D13} as well.

\medskip

\noindent{\bf Conjecture 1:} If $R$ is a periodic semi-perfect ring, then the matrix ring $\mathbb{M}_n(R)$ is periodic.

\medskip

An other claim of interest, which is closely related to the results obtained so far, is the following one:

\begin{prop}\label{sumpotents}
\begin{enumerate}
\item[\rm{(1)}] If $R$ is a potent ring, then every matrix in $\mathbb{M}_n(R)$ is the sum of a tripotent and a potent for all $n\in \mathbb{N}$.
\item[\rm{(2)}] If $R$ is a potent ring and $3\in U(R)$, then every matrix in $\mathbb{M}_n(R)$ is the sum of an idempotent and a potent for all $n\in \mathbb{N}$.
\end{enumerate}
\end{prop}

\begin{proof} Let $n$ be an arbitrary natural number and take $A\in \mathbb{M}_{n}(R)$. Consider now the subring $S$ of $R$ generated by the elements of the matrix $A$. One straightforwardly verifies that the ring $S$ is necessarily finite. Hence, one writes that, $S\cong P_{1} \times \cdots \times P_{m}$ for some finite fields $P_{i}$ with $1\leq i \leq m$; $m\in \mathbb{N}$. Thus, the two assertions of the theorem follow directly from combination of \cite[Theorem 1]{laa-2021} and \cite[Theorem 14]{smz-2021}, as stated.
\end{proof}

So, a non-trivial question which immediately arises is the following (compare also with Conjecture 1 posed above): Suppose that $R$ is an arbitrary periodic ring. Is then, for each $n\geq 1$, any matrix in the ring $\mathbb{M}_n(R)$ a sum of two potent elements? In this direction, Proposition~\ref{sumpotents}(1) could eventually be extended to weakly $2$-primal periodic rings: in fact, consulting with Proposition~\ref{proposition2.1}, the quotient-ring $R/J(R)$ is potent and the ideal $J(R)$ is locally nilpotent. Therefore, each matrix in the factor-ring $\mathbb{M}_n(R)/\mathbb{M}_n(J(R))\cong \mathbb{M}_n(R/J(R))$ is a sum of a tripotent and a potent for any size $n$. Moreover, referring to \cite{BM}, the ideal $J(\mathbb{M}_n(R))=\mathbb{M}_n(J(R))$ is nil. So, what only remains to show is that such a decomposition can be preserved in the whole ring $\mathbb{M}_n(R)$, but the procedure of lifting potent elements and even tripotents is rather more complicated than the standard trick for lifting idempotents (compare with the methods described in \cite{DGL} and \cite{khurana2021lifting}). At this stage, we are unable to do that, so that this should be treated as a left-open problem.\\

Let $A$, $B$ be two rings and $M$, $N$ be the ($A,B$)-bimodule and ($B,A$)-bimodule, respectively. Also, we consider the bilinear maps $\phi:M \otimes_B N \rightarrow A$ and $\psi: N \otimes_A M \rightarrow B$ that apply to the following properties
\[Id_M \otimes_B \psi = \phi \otimes_A Id_M, \quad  \quad Id_N \otimes_A \phi=\psi \otimes_B Id_N. \]

For $m \in M$ and $n \in N$, we define $mn:=\phi(m \otimes n)$ and $nm:= \psi(n \otimes m)$.

Thus, the $4$-tuple $R=\begin{pmatrix} A & M \\ N & B \end{pmatrix}$ becomes to an associative ring equipped with the obvious matrix operations, which is called a \emph{Morita context ring}. Denote the two-sided ideals $Im\phi$ and $Im\psi$ to $MN$ and $NM$, respectively, that are called the \emph{trace} ideals of the Morita context.

\medskip

We are now in a position to establish the following statement, which two different proofs are rather more easy and conceptual than \cite[Theorem 10.1.16]{CS}.

\begin{thm}\label{theorem2.9}
Let
$R =\begin{pmatrix}
  A & M\\
  N & B
\end{pmatrix}$ be a Morita context ring such that $MN$ and $NM$ are nilpotent ideals of $A$ and $B$, respectively. Then, $R$ is periodic if and only if both $A$ and $B$ are periodic.
\end{thm}

\begin{proof} \textbf{Method 1:}

Let $A=
  \left(
  \begin{array}{cc}
  a & m\\
  n & b
  \end{array}
  \right) \in R$.
Also, let $
  \overline{a}=
  \left(
  \begin{array}{cc}
  a & 0\\
  0 & 0
  \end{array}
  \right),
$ $
  \overline{m}=
  \left(
  \begin{array}{cc}
  0 & m\\
  0 & 0
  \end{array}
  \right)
$,
$
  \overline{n}=
  \left(
  \begin{array}{cc}
  0 & 0\\
  n & 0
  \end{array}
  \right)
$
and
$
   \overline{b}=
  \left(
  \begin{array}{cc}
  0 & 0\\
  0 & b
  \end{array}
  \right).
$
Since $MN$ and $NM$ are nilpotent ideals, and rings
$A$ and $B$ are periodic, then the following sets $\{\overline{a}^{t} \mid t\in\mathbb{N} \},$ $\{\overline{b}^{t} \mid t\in\mathbb{N} \},$
$\{\overline{a}^{t_1}\overline{m}\overline{b}^{t_2}\overline{n}\ldots \overline{m}\overline{b}^{t_k}          \mid k\in\mathbb{N}, t_1,\ldots, t_k\in \mathbb{N}_0 \}$,
$\{\overline{a}^{t_1}\overline{m}\overline{b}^{t_2}\overline{n}\ldots \overline{n}\overline{a}^{t_k}          \mid k\in\mathbb{N}, t_1,\ldots, t_k\in \mathbb{N}_0 \}$,
$\{\overline{b}^{t_1}\overline{n}\overline{a}^{t_2}\overline{m}\ldots \overline{n}\overline{a}^{t_k}          \mid k\in\mathbb{N}, t_1,\ldots, t_k\in \mathbb{N}_0 \}$ and $\{  \overline{b}^{t_1}\overline{n}\overline{a}^{t_2}\overline{m}\ldots \overline{m}\overline{b}^{t_k}          \mid k\in\mathbb{N}, t_1,\ldots, t_k\in \mathbb{N}_0 \}$ are finite.
Thus, the set $\{A^{t} \mid t\in\mathbb{N} \}$ is also finite and, therefore, simple calculations lead us to the inspection that, for some different $m,n \in \mathbb{N}$, the equality $A^n=A^m$ holds.

\textbf{Method 2:}

($\Rightarrow$). As both $A$ and $B$ are subrings of $R$, they are obviously periodic.

($\Leftarrow$). If $A$ and $B$ are both periodic, then they are strongly $\pi$-regular. Utilizing \cite[Theorem 3.5]{TLZ}, the ring $R$ is too strongly $\pi$-regular. So, $J(R)$ is nil. But, using \cite[Lemma 3.1]{TLZ}, we have that $R/J(R) \cong A/J(A) \times B/J(B)$. Hence, our assumption and \cite[Lemma 2.12]{BCL} imply that $R/J(R)$ is periodic. Finally, \cite[Proposition 2.4]{BCL} gives the desired result.
\end{proof}

Let $R$, $S$ be two rings, and let $M$ be an $(R,S)$-bi-module such that the operation $(rm)s=r(ms)$ is valid for all $r\in R$, $m\in M$ and $s\in S$. Given such a bi-module $M$, we can set

\[T(R,S,M)=\begin{pmatrix}
R & M \\
0 & S
\end{pmatrix}=\left\{ \begin{pmatrix}
                        r & m \\
                        0 & s
                       \end{pmatrix} : r \in R, m \in M, s \in S  \right\},
\]

\noindent where it forms a ring with the usual matrix operations. The so-stated formal matrix $T(R,S,M)$ is called a {\it formal triangular matrix ring}. In Theorem \ref{theorem2.9}, if we set $N=\{0\}$, then we will obtain the following corollary.

\begin{cor}[{\cite[Theorem 2.13]{BCL}}]
Let $R$, $S$ be two rings and let $M$ be an $(R,S)$-bi-module. Then, the formal triangular matrix ring $T(R,S,M)$ is periodic if and only if $R$ and $S$ are periodic rings.
\end{cor}

By a simple induction, this corollary can be extended to the more general situation of generalized triangular matrix rings. Such rings are denoted by $T(R_i,M_{ij})| 1 \leq i \leq j \leq n)$, where $R_i$ and $M_{ij}$ are, respectively, periodic rings and $(R_i,R_j)$-bi-modules equipped with maps guaranteeing that the multiplication of the matrices is well-defined and satisfies the standard associative property. If, however, $n=3$, this gives that the generalized triangular matrix ring
$$S=\begin{pmatrix}
R_1 & M_{12} & M_{13} \\
0   & R_2    & M_{23} \\
0   & 0      & R_3
\end{pmatrix}$$
is periodic, because
$S=\begin{pmatrix}
A & M \\
0 & R_3
\end{pmatrix}$ with
$A=\begin{pmatrix}
R_1 & M_{12} \\
0 & R_2
\end{pmatrix}$ and
$M=\begin{pmatrix}
M_{13} \\
M_{23}
\end{pmatrix}$, where $R_1$, $R_2$, and $R_3$ are periodic, whereas $M_{12}$, $M_{23}$, $M_{13}$ are, respectively, $(R_1,R_2)$-, $(R_2,R_3)$, $(R_1,R_3)$-bi-modules equipped with the map $\varphi: M_{12} \times M_{23} \rightarrow M_{13}$.

As usual, the notation $\mathbb{T}_n(R)$ stands for the (upper) triangular $n\times n$ matrix ring over an arbitrary ring $R$. Now, we have the following easy corollary.

\begin{cor}
Let $R$ be a ring and $n\geq 1$ is a natural number. Then, $\mathbb{T}_n(R)$ is periodic if and only if $R$ is  periodic.
\end{cor}

Given a ring $R$ and a central element $s$ of $R$, the $4$-tuple $\begin{pmatrix} R & R \\ R & R \end{pmatrix}$ becomes a ring with addition defined componentwise and with multiplication defined by
\[\begin{pmatrix} a_1 & x_1 \\ y_1 & b_1 \end{pmatrix} \begin{pmatrix} a_2 & x_2 \\ y_2 & b_2 \end{pmatrix} =
\begin{pmatrix} a_1a_2+sx_1y_2 & a_1x_2+x_1b_2 \\ y_1a_2+b_1y_2 & sy_1x_2+b_1b_2 \end{pmatrix}. \]
This ring is denoted by $K_s(R)$. A Morita context $\begin{pmatrix} A & M \\ N & B \end{pmatrix}$ with $A=B=M=N=R$ is called a {\it generalized matrix ring} over $R$. It was observed in \cite{KrylovGMR} that a ring $S$ is a generalized matrix ring over $R$ if and only if $S=K_s(R)$ for some $s\in \mathrm{C}(R)$. Here $MN=NM=sR$, so that $MN \subseteq J(A) \Leftrightarrow s \in J(R)$, $NM \subseteq J(B) \Leftrightarrow s \in J(R)$, and $MN,NM$ are nilpotent $\Leftrightarrow s$ is a nilpotent. Thus, Theorem \ref{theorem2.9} has the following consequence, too.

\begin{cor}\label{corollary2.10}
Let $R$ be a ring, and $s \in C(R) \cap \mathrm{Nil}(R)$. Then the formal matrix ring $K_s(R)$ is periodic if and only if $R$ is periodic.
\end{cor}

For $n \geq 2$ and for $s\in C(R)$, the $n \times n$ formal matrix ring over $R$, associated with $s$ and denoted by $\mathbb{M}_n(R;s)$ can be defined same as in \cite{TZ}. Now, we have the following result, which proof is quite more transparent than \cite[Corollary 10.1.17]{CS}.

\begin{thm}\label{theorem2.11}
Let $R$ be a ring with $s\in C(R)\cap \mathrm{Nil}(R)$ and $n\geq 2$. Then $\mathbb{M}_n(R;s)$ is periodic if and only if $R$ is periodic.
\end{thm}

\begin{proof}
We shall use induction on $n$. If $n=2$, then $\mathbb{M}_2(R;s)=K_{s^2}(R)$. So, the claim is true in view of Corollary \ref{corollary2.10}. Suppose now $n >2$ and assume that the claim holds for $\mathbb{M}_{n-1}(R;s)$. Letting $A=\mathbb{M}_{n-1}(R;s)$, we observe that $\mathbb{M}_n(R;s)=\begin{pmatrix}
A & M \\
N & R
\end{pmatrix}$ is a Morita context, where

\[M=\begin{pmatrix}
M_{1n} \\
M_{2n}  \\
\vdots  \\
 M_{n-1,n} \\
\end{pmatrix} , \quad N=\begin{pmatrix}
 M_{n1} M_{n2} \ldots M_{n,n-1}
\end{pmatrix} , \quad M_{in}=M_{ni}=R \quad 1 \leq i \leq n-1. \]
Moreover, for $x=\begin{pmatrix}
                   x_{1n} \\ x_{2n} \\ \vdots \\ x_{n-1,n}
                 \end{pmatrix}$ and $y=\begin{pmatrix}
                                         y_{n1} y_{n2} \ldots y_{n,n-1}
                                       \end{pmatrix}$, we calculate that

\begin{equation} \label{star1}
xy=\begin{pmatrix}
s^2x_{1n}y_{n1} & sx_{1n}y_{n2}    & \ldots & sx_{1n}y_{n,n-1} \\
sx_{2n}y_{n1}    & s^2x_{2n}y_{n2} & \ldots & sx_{2n}y_{n,n-1} \\
        \vdots         &       \vdots           & \ldots & \vdots \\
sx_{n-1,n}y_{n1} & sx_{n-1,n}y_{n2} & \ldots & s^2x_{n-1,n}y_{n,n-1}\\
\end{pmatrix} \in A,
\end{equation}

\begin{equation} \label{star2}
yx=\begin{pmatrix}
y_{n1} y_{n2} \ldots y_{n,n-1}
\end{pmatrix} \begin{pmatrix}
 x_{1n} \\ x_{2n} \\ \vdots \\ x_{n-1,n} \\
\end{pmatrix}=\sum_{k=1}^{n-1} s^2y_{nk}x_{kn} \in R.
\end{equation}

But by (\ref{star1}) we have $MN \subseteq sA$, and by (\ref{star2}) we have $NM=s^2R$. It is, however, obvious that $MN$ and $NM$ are nilpotent ideals of $A$ and $R$, respectively. Now, the claim is sustained from the induction hypothesis combining it with Theorem \ref{theorem2.9}.
\end{proof}

The following technicality is useful for proving our main result stated below.

\begin{lem}\label{lemma2.1} 
Let $R$ be an Abelian periodic ring. If $J(\mathbb{M}_n(R))$ is a nil-ideal of $\mathbb{M}_n(R)$ for some natural number $n$, then the ring $\mathbb{M}_n(R)$ is periodic.
\end{lem}

\begin{proof} 
\textbf{Method 1:} Let $f:R\to R/J(R)$ be a natural homomorphism and $\phi: \mathbb{M}_n(R)\to \mathbb{M}_n(R/J(R))$ the corresponding ring homomorphism induced by $f$. Given $A\in \mathbb{M}_n(R)$ and consider the subring $S$ of the ring $R$, generated by the components of the matrix $A$. Since $R/J(R)$ is a potent ring, one sees that $f(S)$ is a finite potent subring of $R/J(R)$. Therefore, $f(S)= e_1f(S)\oplus\ldots \oplus e_kf(S)$, where $e_1,\ldots, e_k\in f(S)$ is a family of pairwise orthogonal idempotents and $e_if (S)$ is a finite field. Moreover, there exist pairwise orthogonal idempotents, say $f_1,\ldots, f_n\in S$, such that $S= f_1S\oplus\ldots \oplus f_nS$ and $\phi(f_i)=e_i$ for each $1\leq i \leq n$. Thus $A=A_1+\ldots +A_k,$ where $A_i\in \mathbb{M}_n(f_iS)$. Let us now we set
\[\begin{split}
q=&[[\mid e_{1}f(S)\mid-1, \ldots, \mid e_{1}f(S)\mid^n-1],\ldots, \\
  &[[\mid e_{k}f(S)\mid-1, \ldots, \mid e_{k}f(S)\mid^n-1]]]+1.
\end{split}\]
Since $J(R)=\mathrm{Nil}(R)$, it is easy to see that, for each $1\leq i \leq n$, the equality $\mathrm{Nil}(f_iS)=J(f_iS)$ holds. Consequently, $$J(\mathbb{M}_n(f_iS))=\mathbb{M}_n(\mathrm{Nil}(f_iS))\subseteq \mathbb{M}_n(J(R))$$
and $J(\mathbb{M}_n(f_iS))$ is a nil-ideal. So, applying \cite[Theorem 8]{AT}, for each $1\leq i \leq n$ we may deduce that $A_i^{q}-A_i\in \mathrm{Nil}(\mathbb{M}_n(f_iS))$. Hence, for every $A\in \mathbb{M}_n(R)$, we obtain that $A^{q}-A\in \mathrm{Nil}(\mathbb{M}_n(R))$ and thus, owing to \cite[Theorem 3.4]{CD}, we can conclude that $\mathbb{M}_n(R)$ is a periodic ring, as wanted.

\textbf{Method 2:} 
If R is an Abelian periodic ring then $J(R)=Nil(R)$ by \cite[Lemma 5]{Badawi}. Now, the isomorphism $\mathbb{M}_n(R)/\mathbb{M}_n(Nil(R)) \cong \mathbb{M}_n(R/Nil(R))$, Lemma \ref{lemma3.1}, the hypothesis imply that $\mathbb{M}_n(R)$ is periodic.
\end{proof}

We now offer the following statement which, somewhat, is in a parallel to \cite{M}.

\begin{thm}\label{theorem2.7} The following points are equivalent:

\begin{enumerate}
\item[\rm{(1)}] If $R$ is a periodic abelian ring, then $\mathbb{M}_2(R)$ is a periodic ring.

\item[\rm{(2)}] If $R$ is a periodic local ring, then $\mathbb{M}_2(R)$ is a periodic ring.

\item[\rm{(3)}] The K\"othe's Conjecture has a positive solution for all algebras over finite fields.
\end{enumerate}
\end{thm}

\begin{proof} The equivalence $(1) \Leftrightarrow (2)$ follows directly from proof of Lemma \ref{lemma2.1}.

We shall now deal with the implication $(2) \Rightarrow (3)$. Take a nil-algebra $A$ over a finite field $\mathbb{F}$. Write $A^*$ for an algebra over $\mathbb{F}$, which is the Dorroh overring of $A$ (see, e.g., \cite{D13}). In view \cite[Theorem 3.4]{CD}, $A^*$ is a periodic ring. In the presence of condition (2), the algebra $M_2(A^*)$ is a periodic ring. Hence $J(\mathbb{M}_2(A^*))=\mathbb{M}_2(J(A^*))=\mathbb{M}_2(A)$ is a nil algebra.

We now consider the reverse implication $(3) \Rightarrow (2)$. Take a periodic local ring $R$. Then, $J(R)$ is nil and the quotient $R/J(R)$ is a finite field. Letting char$(R/J(R))=p$, then condition (3) implies that $\mathbb{M}_2(J(R/pR))=\mathbb{M}_2(J(R)/pR))$ is nil. Since $(pR)^n=0$ for some $n\in \mathbb{N}$, one infers that $J(\mathbb{M}_2(R))=\mathbb{M}_2(J(R))$ is nil too. Now, Lemma \ref{lemma2.1} yields that the desired result.
\end{proof}

The following comments could be useful.

\begin{rem}\label{Dis} Surprisingly, in accordance with \cite[Theorem 3.1]{M} and point (3) of the last theorem of ours, the famous Diesl's question from \cite{D13} whether or not {\it $R$ being nil-clean implies that $\mathbb{M}_n(R)$ is nil-clean for all $n\geq 1$} is paralleling to the corresponding implication for (abelian, local) periodic rings. On the other side, concerning the final successful resolving of the question of Diesl, we also suggest the following tactic: Let $R$ be a (local, abelian) periodic ring. If we can prove that $J(R)$ is locally nilpotent, then we know that $J(\mathbb{M}_n(R))$ is nil and hence $\mathbb{M}_n(R)$ is periodic by our Lemma~\ref{lemma2.1}. Next, our Theorem~\ref{theorem2.7} combined with \cite[Theorem 3.1]{M} apply to get the desired solution after all.
\end{rem}

\section{Endomorphism Rings}\label{SectionEndomorphismRings}

\medskip

In his seminal paper \cite{FR}, Fuchs along with Rangaswamy do {\it not} succeeded to give a full characterization when the endomorphism ring $\mathrm{E}(G)$ of an arbitrary Abelian group $G$ is $\pi$-regular, excepting the torsion case (see, for more account, \cite{F} and \cite{krylov2013endomorphism}, too). In fact, they actually established some deep results in the general mixed case as well as in the case of torsion-free groups, and a complete description when $G$ is a torsion group by proving that it has to be specially bounded. It is not known to the authors whether or not the (stronger) property of being strongly $\pi$-regular is already studied for endomorphism rings of Abelian groups. For the more general case of clean endomorphism rings of abelian groups we refer to \cite{GV}. In this vein, we just mention that it was showed in \cite[Proposition 10]{BCDM} that the endomorphism ring of an arbitrary Abelian group is strongly nil-clean exactly when the group is a cyclic $2$-group as well as in \cite[Proposition 9]{BCDM} that the endomorphism ring of a finite rank Abelian group is nil-clean exactly when it is a finite $2$-group.

We begin our work here with the following two statements, which show some relationship between the results from the preceding section and these from the present one.

\begin{thm} For a ring $R$, the following two conditions are equivalent:
\begin{enumerate}
 \item[\rm{(1)}] For each $n\in \mathbb{N}$, the ring $\mathbb{M}_n(R)$ is periodic;

 \item[\rm{(2)}] The endomorphism ring $\mathrm{E}(M)$ of every finitely generated right $R$-module $M$ is a periodic ring.
\end{enumerate}
\end{thm}

\begin{proof}
The implication $(2)\Rightarrow (1)$ is rather obvious, so we leave it to be checked by the interested reader.

$(1)\Rightarrow (2)$. Let $M$ be a finitely generated right $R$-module. Then, for some finitely generated free right $R$-module $F$, there clearly exists an epimorphism (i.e., a surjective homomorphism) $f: F\to M$. Moreover, if $\phi\in \mathrm{E}(M)$, then for some $\phi'\in \mathrm{E}(F)$ the equality $f\phi'=\phi f$ holds. Thus, according to our assumption, for some $m,n \in \mathbb{N}$ the equality $\phi'^m=\phi'^n$ is valid. Therefore, one sees that the equalities
$$\phi ^nf= f\phi'^n=f\phi'^m=\phi^mf$$
\noindent are true. However, since $f$ is an epimorphism, one concludes that $\phi^n=\phi^m$, as needed.
\end{proof}

The following consequence follows at once from the previous theorem and Theorem \ref{theorem2.5} (see also Corollary~\ref{corollary2.6}).

\begin{cor}
If $R$ is a right (resp., left) perfect periodic ring (in particular, a periodic Artinian ring) and $M$ is a finitely generated right $R$-module, then the ring $\mathrm{E}(M)$ is periodic.
\end{cor}

The next technical statement is pivotal.

\begin{lem}
Let $V$ be a vector space over a field $F$. Then, $\mathrm{E}(V)$ is a periodic ring if and only if $V$ is a finite-dimensional vector space and $F$ is an algebraic extension of a finite field.
\end{lem}

\begin{proof}
Let us show in a way of contradiction that if $V$ is an infinite-dimensional vector space, then $\mathrm{E}(V)$ is not a periodic ring. Obviously, it suffices to consider the case when $V$ is a countable-dimensional space. Let $e_1, e_2, \ldots, e_k$ be a basis of the vector space $V$. Consider now a linear operator $f:V\to V$ that acts on basis vectors according to the equalities $f(e_1)=0, f(e_2)=e_1, f(e_3)=e_2, \ldots , f(e_k)=e_{k-1}$. It is, therefore, an easy technical exercise to see that $f^n\neq f ^m$ for every natural numbers $m\neq n$, as required.

To show now the second part, namely that $F$ is an algebraic extension of a finite field, we firstly observe that the field $F$ is a periodic ring. So, it is obvious that the characteristic of this field is non-zero and is equal to some prime number, say $p$.  Furthermore, since each element of the field $F$ is periodic, the field $F$ is really an algebraic extension of the finite simple field
$\mathbb{F}_p$, as claimed.
\end{proof}

Since for bounded Abelian $p$-groups $G$, the ring $\mathrm{E}(G)/p\mathrm{E}(G)$ is always isomorphic to the ring of row-finite matrices over the simple field $\mathbb{F}_p$, then the following assertion follows directly from the previous lemma.

\begin{lem} \label{lemma3.10}
For any bounded Abelian $p$-group $G$, the ring $\mathrm{E}(G)$ is periodic if and only if $G$ is a finite group.
\end{lem}

Our main result here is as follows. We prove it in two methods.

\begin{thm} \label{theorem3.11}
For any Abelian group $G$, the ring $\mathrm{E}(G)$ is periodic if and only if $G$ is a finite group.
\end{thm}

\begin{proof} \textbf{Method 1:}

($\Rightarrow$). Let $\mathrm{E}(G)$ be periodic. Take the endomorphism $\varphi:G \rightarrow G$ defined by $\varphi(x)=2x$. As $\mathrm{E}(G)$ is periodic, there will exist two natural numbers, say $n$ and $m$, such that ${\varphi}^n={\varphi}^m$ ($n>m$). Thus $(2^n-2^m)x=0$ for any $x \in G$, i.e., $G$ is of finite exponent. However, by \cite[Theorem 8.4]{fuchs}, $G$ is a finite direct sum of bounded $p_i$-groups $G_{p_i}$ for some different primes $p_i$. So, the isomorphism $\mathrm{E}(G) \cong \prod_{i=1}^n \mathrm{E}(G_{p_i})$ yields that $\mathrm{E}(G_{p_i})$ is periodic. Finally, referring to Lemma \ref{lemma3.10}, $G_{p_i}$ is a finite group. Therefore, $G$ is a finite group, as asserted.

\medskip

\textbf{Method 2:}

($\Rightarrow$). Let $\mathrm{E}(G)$ be periodic. Similarly to the first method demonstrated above, $G$ is of finite exponent. By \cite[Proposition 2.5]{BCL}, the ring $\mathrm{E}(G)$ is Dedekind finite (or, in other terms, directly finite). Now, \cite[Exercise 15, p.121]{krylov2013endomorphism} implies $G$ is finite.

($\Leftarrow$). If $G$ is finite, then $\mathrm{E}(G)$ is finite. Now, by Lemma \ref{lemma3.1}, the ring $\mathrm{E}(G)$ is periodic, as stated.
\end{proof}

Indenting to expand the aforementioned results from \cite{BCDM} concerning the strong nil-cleanness of $\mathrm{E}(G)$, we will now consider the more large class of strongly $m$-nil-clean rings for some integer $m>1$ that are rings for which every element is a sum of an $m$-potent and a nilpotent which commute each other. Obviously, strongly $m$-nil clean rings are always periodic. In particular, if $m=2$, these rings are just the strongly nil-clean ones.

Let $R$ be a ring. If $r\in R$ is nilpotent, we denote its order of nilpotency by $\mathrm{In}(r) = \mathrm{min}\{k\mid r^k = 0\}$, and if $r$ is not nilpotent we just put $\mathrm{In}(r)=0$. Thereby, the degree of nilpotency $\mathrm{In}(R)$ of $R$ is defined by $\mathrm{In}(R) = \mathrm{sup}_{r\in R}~\mathrm{In}(r)$.

\begin{thm}\label{theorem3.12}
Let $m>1$ be a natural number and $G$ an Abelian group. Then, the following statements are equivalent:

\begin{enumerate}
\item[\rm{(1)}] $\mathrm{E}(G)$ is a strongly $m$-nil clean ring;

\item[\rm{(2)}] $G$ is a finite group and, for each $p$-primary component $G_p$ of the group $G$, it must be that $(p^i-1) \mid (m-1)$ for every $1\leq i\leq \mathrm{In}(\mathrm{E}(G_p)/J(\mathrm{E}(G_p)))$.

\end{enumerate}
\end{thm}

\begin{proof} Since every strongly $m$-nil clean ring is periodic, in conjunction with Theorem \ref{theorem3.11} we may assume without loss of generality in proving the equivalence $(1) \Leftrightarrow (2)$ that the group $G$ is finite.

And so, suppose that $G_p$ is the $p$-primary component of $G.$ By \cite{fuchs}, one may write that $G_p\cong \mathbb{Z}(p^{k_1})^{(n_1)}\oplus\ldots \oplus \mathbb{Z}(p^{k_t})^{(n_t)}$, where $k_1<\ldots<k_t$. It is now easy to detect that $$\mathrm{E}(G_p)/J(\mathrm{E}(G_p))\cong \mathbb{M}_{n_1}(\mathbb{F}_p)\times\ldots \times  \mathbb{M}_{n_t}(\mathbb{F}_p)$$
and that $J(\mathrm{E}(G_p))$ is a nilpotent ideal. Therefore, it follows from \cite[Theorem 8]{AT} that $\mathrm{E}(G_p)$ is strongly $m$-nil clean ring if and only if $(p^i-1)\mid (m-1)$ for every $1\leq i\leq \mathrm{max}\{n_1,\ldots, n_t\}=\mathrm{In}(\mathrm{E}(G_p)/J(\mathrm{E}(G_p))$. Thus, the wanted equivalence follows directly from the fact that the isomorphism $\mathrm{E}(G)\cong \prod_p \mathrm{E}(G_p)$ is fulfilled.
\end{proof}

Some partial cases of the last result are also of interest. Although the following statement follows automatically from the proof of the previous theorem and \cite[Corollary 19]{AT}, we prefer to prove the following result directly by illustrating a more transparent and conceptual proof.

\begin{cor}\label{corollary3.13}
Let $m$ be a positive integer with $m \not\equiv 1$ (mod 3), $m \not\equiv 1$ (mod 8). Then the ring $\mathrm{E}(G)$ is strongly $m$-nil clean if and only if $G\cong \oplus_{i=1}^n \oplus_{j=1}^{t_i} \mathbb{Z}_{{p_i}^{k_j}}$, where $p_1,p_2, \cdots,p_n$ are prime factor of $|G|$ and $(p_i-1)|(m-1)$ for any $1 \leq i \leq n$.
\end{cor}

\begin{proof} ($\Rightarrow$). Suppose that $\mathrm{E}(G)$ is strongly $m$-nil clean. So, Theorem \ref{theorem3.11} implies that $G$ is a finite group. Hence, $G=\oplus_{i=1}^{n} G_{p_i}$, where $G_{p_i} \cong \mathbb{Z}({p_i}^{k_1})^{(n_1)}\oplus\ldots \oplus \mathbb{Z}({p_i}^{k_{t_i}})^{(n_{t_i})}$ and $k_1<\ldots<k_{t_i}$. It is clear that $n$ is the number of prime factors of $|G|$. But $$\mathrm{E}(G_{p_i})/J(\mathrm{E}(G_{p_i}))\cong \mathbb{M}_{n_1}(\mathbb{F}_{p_i})\times\ldots \times  \mathbb{M}_{n_{t_i}}(\mathbb{F}_{p_i}).$$
The hypothesis yields that $\mathbb{M}_{n_1}(\mathbb{F}_{p_i}), \ldots ,  \mathbb{M}_{n_{t_i}}(\mathbb{F}_{p_i})$ are strongly $m$-nil clean. Thus, it follows from \cite[corollary 19]{AT} that $n_1=n_2= \ldots =n_{t_i}=1$. Therefore, $G \cong \oplus_{i=1}^n \oplus_{j=1}^{t_i} \mathbb{Z}_{{p_i}^{k_j}}$.

($\Leftarrow$). Suppose that $G \cong \oplus_{i=1}^n \oplus_{j=1}^{t_i} \mathbb{Z}_{{p_i}^{k_j}}$, where $(p_i-1)|(m-1)$ for any $1 \leq i \leq n$. So, $\mathrm{E}(G) \cong \prod_{i=1}^n \prod_{j=1}^{t_i} \mathrm{E}(\mathbb{Z}_{{p_i}^{k_j}})$. As the isomorphism $\mathrm{E}(\mathbb{Z}_{{p_i}^{k_j}}) \cong \mathbb{Z}_{{p_i}^{k_j}}$ is valid, one detects that the ring $$\frac{\mathrm{E}(\mathbb{Z}_{{p_i}^{k_j}})}{p_i\mathrm{E}(\mathbb{Z}_{{p_i}^{k_j}})} \cong \frac{\mathbb{Z}_{{p_i}^{k_j}}}{p_i\mathbb{Z}_{{p_i}^{k_j}}} \cong \mathbb{F}_{p_i}$$
is strongly $m$-nil clean. By considering the denominator, which is obviously a nilpotent ideal, \cite[Theorem 1.11]{BMA} guarantees that $\mathrm{E}(\mathbb{Z}_{{p_i}^{k_j}})$ is a strongly $m$-nil clean ring. Consequently, $\mathrm{E}(G)$ is too strongly $m$-nil clean, as promised.
\end{proof}

\begin{cor}
For an even integer $m$ with $m \not\equiv 1$ (mod 3), the ring $\mathrm{E}(G)$ is strongly $m$-nil clean if and only if $G\cong \oplus_{j=1}^t \mathbb{Z}_{2^j}$.
\end{cor}

\begin{proof}
Taking into account that $2 \in \mathrm{Nil}(\mathrm{E}(G))$, the group $G$ has to be a finite $2$-group. Furthermore, Corollary \ref{corollary3.13} gives the wanted result.
\end{proof}

\section{Tensor Product of Algebras}\label{SectionTensorProductOfAlgebra}

\medskip

Throughout this section, we shall assume that $k$ is a commutative ring. We begin our considerations with the following technical convention.

\begin{lem}\label{lemma5.1}
Suppose that $A$ and $B$ are two $k$-algebras such that, as rings, they are periodic commutative rings. Then the tensor product $k$-algebra $A \otimes_k B$ is periodic.
\end{lem}

\begin{proof}
It is sufficient to show that $(A \otimes_k B)/Nil(A \otimes_k B)$ is a periodic $k$-algebra. Suppose that $\overline{\sum_{i=1}^t x_i \otimes y_i}$ is an arbitrary element of $(A \otimes_k B)/\mathrm{Nil}(A \otimes_k B)$. Since $A$ and $B$ are periodic, $x_i=e_i+a_i$, where  ${e_i}^{n_i}=e_i$ for some $n_i \in \mathbb{N}$, $a_i \in \mathrm{Nil}(A)$ and $y_i=f_i+b_i$, where ${f_i}^{m_i}=f_i$ for some $m_i \in \mathbb{N}$, $b_i \in \mathrm{Nil}(B)$ for $1 \leq i \leq t$. It is easy to verify that $\overline{\sum_{i=1}^t x_i \otimes y_i}=\sum_{i=1}^t \overline{e_i \otimes f_i}$. Furthermore, one finds by \cite[Lemma 2.3]{BCL} that ${e_i}^{k(n_i-1)+1}=e_i$, ${f_i}^{k(m_i-1)+1}=f_i$ for any $k \in \mathbb{N}$. If $l_i=(n_i-1)(m_i-1)+1$, then ${e_i}^{l_i}=e_i$, ${f_i}^{l_i}=f_i$. So, ${(\overline{e_i \otimes f_i})}^{l_i}=\overline{{e_i}^{l_i} \otimes {f_i}^{l_i}}=\overline{e_i \otimes f_i}$. Therefore, $\overline{e_i \otimes f_i}$ is a potent element. But
\[(\overline{e_i \otimes f_i})(\overline{e_j \otimes f_j})=\overline{e_ie_j \otimes f_if_j}=\overline{e_je_i \otimes f_jf_i}=(\overline{e_j \otimes f_j})(\overline{e_i \otimes f_i}). \]
Thus, in conjunction with \cite[Lemma 2.6]{BCL}, the element $\overline{\sum_{i=1}^t x_i \otimes y_i}=\sum_{i=1}^t \overline{e_i \otimes f_i}$ is periodic. Consequently, $(A \otimes_k B)/Nil(A \otimes_k B)$ is periodic as a ring, as required.
\end{proof}

We thus obtain the following.

\begin{thm}\label{theorem5.2}
Suppose that $A$ and $B$ are two $k$-algebras such that, as rings, they are weakly $2$-primal periodic rings. Then, $A \otimes_k B$ is a periodic $k$-algebra.
\end{thm}

\begin{proof}
The given assumption implies that $A/\mathrm{Nil}(A)$ and $B/\mathrm{Nil}(B)$ are two potent $k$-algebras. Now, Lemma \ref{lemma5.1} allows us to deduce that $(A/\mathrm{Nil}(A)) \otimes_k (B/\mathrm{Nil}(B))$ is a periodic $k$-algebra. The rest of our argument is quite similar to the one in \cite[Theorem 2.3]{BMA}.
\end{proof}

Our next technical convention is the following.

\begin{lem}\label{lemma5.3}
Suppose that $A$ and $B$ are two $k$-algebras such that they are locally finite as rings. Then, the tensor product $k$-algebra $A \otimes_k B$ is locally finite.
\end{lem}

\begin{proof}
Suppose $S=\{ \sum_{i=1}^{n_1} a_{i1} \otimes b_{i1}, \sum_{i=1}^{n_2} a_{i2} \otimes b_{i2}, \cdots, \sum_{i=1}^{n_t} a_{it} \otimes b_{it} \}$ is a finite set of $A \otimes_k B$. Also, suppose $A^{'}$ is the subalgebra generated by the elements $$\{a_{11},a_{21}, \cdots , a_{n_11},a_{12},a_{22}, \cdots,a_{n_22}, \cdots,a_{1t},a_{2t}, \cdots,a_{n_tt}\}$$
and $B^{'}$ is the subalgebra generated by $$\{b_{11},b_{21}, \cdots , b_{n_11},b_{12},b_{22}, \cdots,b_{n_22}, \cdots,b_{1t},b_{2t}, \cdots,b_{n_tt}\}.$$
Since it is obvious that both $A$ and $B$ are locally finite, it is clear that $A^{'}$ and $B^{'}$ must be finite. So, $A^{'} \otimes_k B^{'}$ is finite itself, too. But $S \subseteq A^{'} \otimes_k B^{'}$ and thus the subalgebra generated by $S$ is finite, as well. Therefore, $A \otimes_k B$ is locally finite, as claimed.
\end{proof}

We finish this subsection with the following result.

\begin{thm}\label{theorem5.4}
Suppose that $A$ and $B$ are two $k$-algebras such that they are right (resp., left) perfect and periodic rings. Then, the tensor product $k$-algebra $A \otimes_k B$ is periodic as a ring.
\end{thm}

\begin{proof}
By assumption, $A/J(A)$ and $B/J(B)$ are semi-simple $k$-algebras. Now, according to the well-known Wedderburn's Structure Theorem for algebras (see, for instance, \cite[Theorem 3.5]{pierce1982}), we have
\[\frac{A}{J(A)} \cong \prod_{i=1}^{t} \mathbb{M}_{n_i}(D_i) \text{  ,  } \frac{B}{J(B)} \cong \prod_{j=1}^l \mathbb{M}_{m_j}(D^{'}_j).\]

\noindent for some division rings $D_i$ and $D^{'}_j$, respectively. It now follows from the periodicity of both $A$ and $B$ that $D_i$, $D^{'}_j$ are periodic and, consequently, they are potent fields for any $1 \leq i \leq t$ and $1 \leq j \leq l$. So, $\mathbb{M}_{n_i}(D_i)$ and $\mathbb{M}_{m_j}(D^{'}_j)$ are locally finite $k$-algebras by consulting with \cite[Corollary 2.3]{HKL}. Furthermore, Lemma \ref{lemma5.3} implies that $\mathbb{M}_{n_i}(D_i) \otimes_k \mathbb{M}_{m_j}(D^{'}_j)$ is a locally finite $k$-algebra for any $1 \leq i \leq t$ and $1 \leq j \leq l$. However, one observes by \cite[Lemma 2.12]{BCL} that
\[\frac{A}{J(A)} \otimes_k \frac{B}{J(B)} \cong \prod_{i=1}^t\prod_{j=1}^l \mathbb{M}_{n_i}(D_i) \otimes_k \mathbb{M}_{m_j}(D^{'}_j)\]
is a periodic $k$-algebra. If now $A$ and $B$ are right (resp., left) perfect, then $J(A)$ and $J(B)$ are locally nilpotent in view of \cite[Proposition 23.15]{L} (see also \cite[Exercise 23.1]{L}). In a way of similarity to the proof of \cite[Theorem 2.3]{BMA}, we may infer that $(i \otimes 1)(J(A) \otimes_k B))+(1 \otimes j)(A \otimes_k J(B))$ is a nil-ideal of $A \otimes_k B$. Moreover, with the aid of \cite[Proposition 2.2]{BMA}, we have the following isomorphism
\[\frac{A \otimes_k B}{(i \otimes 1)(J(A) \otimes_k B)+(1 \otimes j)(A \otimes_k J(B))} \cong \frac{A}{J(A)} \otimes_k \frac{B}{J(B)}.\]
Finally, \cite[Corollary 3.6]{BCL} yields the desired result, as expected.
\end{proof}

It is worthwhile noticing that some valuable results on tensor products of (weakly) nil-clean rings were obtained in \cite{S}.

\section{Weakly Periodic Rings}\label{SectionWeaklyPeriodicRings}

\medskip

It is well known from \cite{CD} that if $R$ is a periodic ring then every element $x\in R$ can be written as $x=e+b$, where $e^n=e$ for some natural number $n \geq 2$ and $b$ is a nilpotent element such that $eb=be$. Likewise, it was proved in \cite{Bell} that if each element in $R$ is expressible as the sum of a potent element and a nilpotent element such that $\mathrm{Nil}(R)$ is commutative, then $R$ is periodic. Later, the authors called the rings whose every element can be represented as the sum of a potent element and a nilpotent element {\it weakly periodic}. It is clear that every periodic ring is weakly periodic, but however the converse manifestly fails -- indeed, according to \cite[Examples 3.1,3.2]{Ster}, we can exhibit a nil-clean ring (and hence a weakly periodic ring) which is {\it not} strongly $\pi$-regular and so {\it not} periodic as well, thus answering the corresponding question raised in \cite{Hirano}. 

Further, for a fixed natural number $m\geq 2$, an element $x\in R$ is called an {\it $m$-potent} whenever $x^m=x$ and so, based on \cite{BMA}, a ring $R$ is called {\it $m$-nil clean} if each its element is presentable as the sum of an $m$-potent element and a nilpotent element. Hence, every $m$-nil clean ring is weakly periodic. In particular, every weakly nil-clean ring, as defined in \cite{DMc} and \cite{BDZ}, is 3-nil clean, because the elements $\pm e$ are obviously tripotents whenever $e$ is an idempotent, so that these rings are too weakly periodic (for an alternative, more conceptual, approach to demonstrate this fact, we refer the interested reader to the main results in \cite{Dnew} and \cite{St}, respectively, which prove that each weakly nil-clean ring $R$ is decomposable as $R\cong R_1\times R_2$, where $R_1$ is a nil-clean ring, and $R_2$ is either $(0)$ or $J(R_2)$ is nil such that $R_2/J(R_2)\cong \mathbb{Z}_3$. Thus, we can automatically employ a simple combination of Propositions~\ref{proposition1.6} and \ref{proposition1.7} listed below to get the wanted assertion after all). For an other important study of some rings of this sort, we are referring to \cite{Di}.

On the other side, in \cite{Ye} the concept of {\it semi-clean rings} was introduced. Recall that an element $x$ in a ring $R$ is {\it periodic} if $x^m=x^n$, where $m$ and $n$ are natural numbers and $m \neq n$. So, a ring $R$ is called semi-clean when each its element is expressible as the sum of a periodic element and an unit element. It is easy to see that any weakly periodic ring is semi-clean. However, the converse fails as it will be shown below.

Moreover, as we know, every periodic ring is strongly $\pi$-regular, whence a logical question is to ask whether every weakly periodic ring is $\pi$-regular. But, unfortunately, by using once again \cite[Example 3.2]{Ster}, one can say that a weakly periodic ring is {\it not} necessarily $\pi$-regular as there exists a nil-clean ring which is {\it not} $\pi$-regular, so another challenging question is of whether or not weakly periodic rings are {\it clean} in the sense of \cite{N}. We, however, hope the answer to be negative too, but at this stage we do not have a concrete example yet.

Motivated by the above discussion, we shall try in this short section to present some basic facts about weakly periodic rings, thus eventually stimulating their further exploration in detail. Specifically, as noted above, our main motivating tool is the following claim, which unambiguously illustrates that the class of weakly periodic rings may have rather thin properties than the larger class of semi-clean rings.

\begin{ex} There exists a semi-clean ring which is {\it not} weakly periodic.
\end{ex}

\begin{proof} According to \cite{Ye} (see also \cite{McGov}), the group ring $\mathbb{Z}_{(p)}C_3$ is semi-clean. However, if we assume in a way of contradiction that $\mathbb{Z}_{(p)}C_3$ is weakly periodic, then $\mathbb{Z}_{(p)}$ will be weakly periodic by Theorem \ref{theorem3.10} quoted below. So, $J(\mathbb{Z}_{(p)})=p\mathbb{Z}_{(p)}$ is nil, a contradiction.

Alternatively, if $R$ is a weakly periodic ring, then the power series ring, $R[[x]]$, is semi-clean by \cite[Proposition 3.3]{Ye}, but as we will see in the next example, $R[[x]]$ is not weakly periodic.
\end{proof}

\begin{ex} For any ring $R$, both the Laurent polynomial ring $R[x,x^{-1}]$ and the power series ring $R[[x]]$ are {\it not} weakly periodic. Also, the polynomial ring $R[x]$ is {\it not} periodic.
\end{ex}

\begin{proof} If we assume, by way of contradiction, that the ring $R[x,x^{-1}]$ is weakly periodic, then the isomorphism $R[x,x^{-1}] \cong RG$, where $G$ is an infinite cyclic group, along with Theorem~\ref{theorem3.10} (2), stated and proved below, will imply that $G$ is a torsion group, which is a contradiction to our assumption, thus giving the wanted claim.

Supposing now that $R[[x]]$ is weakly periodic, we know that $J(R[[x]])=\{a+xf(x) ~|~ a \in J(R)~ \mathrm{and} ~f(x) \in R[[x]]\}$ (see, for instance, Exercise 5.6 in \cite{L}), and so it is clear that $x \in J(R[[x]])$. Therefore, $J(R[[x]])$ is not nil, a contradiction, thus getting the desired claim.

Next, if $R[x]$ is periodic, then it is strongly clean by the usage of \cite[Proposition 2.5]{BCL}, and thus $R[x]$ must be clean. But then this contradicts \cite[Remark 2.8]{kanwar2015clean}, and hence $R[x]$ cannot be periodic, as claimed.
\end{proof}

A question which naturally arises is of whether the polynomial ring $R[x]$ is weakly periodic or not. At this stage, we unable to prove or disprove this question.

\medskip

We continue our work in this section with the following useful technicalities which show some fundamental properties of weakly periodic rings (some other background material can be found in \cite{CS}).

\begin{prop}\label{proposition1.6}
Let $R$ be a weakly periodic ring. Then, the following three claims hold:
\begin{enumerate}
\item[\rm{(1)}]
$R$ has positive characteristic.
\item[\rm{(2)}]
If $\mathrm{char}(R)={p_1}^{n_1}{p_2}^{n_2} \cdots {p_k}^{n_k}$, then $R \cong \prod_{i=1}^k R_i$, where $R_i\cong R/({p_i}^{n_i}R)$ for each index $i$.
\item[\rm{(3)}]
A finite direct product of weakly periodic rings is also weakly periodic.
\end{enumerate}
\end{prop}

\begin{proof}
(1) Knowing that $J(R)$ is nil, so it turns out that the characteristic of $R$ is finite and thus positive.

\noindent (2) Since ${p_i}^{n_i}R+{p_j}^{n_j}R=R$ and ${p_i}^{n_i}R \cap {p_j}^{n_j}R=(0)$, with the Chinese Remainder Theorem at hand, we have that $R \cong \prod_{i=1}^k R_i$, where $R_i\cong R/({p_i}^{n_i}R)$ for each index $i$.

\noindent (3) Let $R=\prod_{i=1}^n R_i$, where, for any $1 \leq i \leq n$, the direct component $R_i$ is weakly periodic. Given $r=(r_1,r_2, \cdots r_n) \in R$ with $r_i \in R_i$. So, $r_i=e_i+b_i$, where ${e_i}^{n_i}=e_i$ and ${b_i}^{m_i}=0$ for any $1 \leq i \leq n$. Consulting with \cite[Lemma 2.3(1)]{BCL}, one detects that ${e_i}^{k(n_i-1)+1}=e_i$ for any $k \in \mathbb{N}$. If, however, we put $l=\prod_{j=1}^n (n_j-1)+1$ and $m=max\{m_j\}_{j=1}^n$, then, for any $1 \leq i \leq n$, we inspect that ${e_i}^l=e_i$, $b_i^m=0$. This readily implies that $r=e+b$, where $e=(e_1,e_2, \cdots e_n)$, $b=(b_1,b_2, \cdots b_n)$ and $e^l=e$, $b^m=0$. Therefore, $R$ is weakly periodic, as wanted.
\end{proof}

\begin{prop}\label{proposition1.7}
Let $R$ be a ring and let $I$ be a nil-ideal of $R$. Then, the ring $R$ is weakly periodic if and only if the quotient-ring $R/I$ is weakly periodic.
\end{prop}

\begin{proof}
($\Rightarrow$). It is obvious, so we omit the details.

($\Leftarrow$). Suppose that $a\in R$. By assumption, $\overline{a}=\overline{e}+\overline{b}$, where ${\overline{e}}^m=\overline{e}$ for some $m \in \mathbb{N}$ and $\overline{b}\in \mathrm{Nil}(R/I)$. Utilizing \cite[Corollary 6]{khurana2021lifting}, there will exist a potent element $f\in R$ such that $\overline{f}=\overline{e}$. Now, as $I$ is a nil-ideal, the relation $a-f\in \mathrm{Nil}(R)$ follows. Since $a=f+(a-f)$, $f$ is a potent element and $a-f$ is a nilpotent element, we deduce that $a$ is weakly periodic. Therefore, $R$ is a weakly periodic ring, as stated.
\end{proof}

\medskip

We now shall be concerned with Morita contexts by establishing the following assertion. Before that, the next technical claim is true.

\begin{lem} \label{lemma5.4}
Let $T$ be a ring such that $T=S+K$, where $S$ is a subring of $R$ and $K$ is a nil ideal of $T$. If $S$ is weakly periodic, then $T$ is weakly periodic.
\end{lem}

\begin{proof}
Let $t \in T$. By hypothesis we may write that $t=s+x$, where $s\in S$ and $x\in K$. As $S$ is weakly periodic, one writes that $s=a+b$, where $a^n=a$ for some integer $n>1$ and $b^l=0$ for some integer $l>1$. So, $t=a+b+x$. Set $c=b+x$. Working in $R/K$ and lifting to $R$, we see $c^l \in K$. Since $K$ is nil, we have ${(c^l)}^m=0$ for some $m \geq 1$. Therefore, $T$ is weakly periodic.  
\end{proof}

We now have the validity of the following.

\begin{thm} \label{theorem5.5}
Let
$R =\begin{pmatrix}
  A & M\\
  N & B
  \end{pmatrix}$ be a Morita context ring such that $MN$ and $NM$ are nilpotent ideals of $A$ and $B$, respectively. Then, $R$ is weakly periodic if and only if both $A$ and $B$ are weakly periodic.
\end{thm}

\begin{proof}
($\Rightarrow$). Suppose that $R$ is weakly periodic. It is clear that $MN \subseteq J(A)$ and $NM \subseteq J(B)$. We have
$J(R)=\begin{pmatrix}
          J(A) & M \\
          N    & J(B)
       \end{pmatrix}$ and $R/J(R) \cong A/J(A) \times B/J(B)$ by the utilization of \cite[Lemma 3.1]{TLZ}.
Since $J(R)$ is nil, it is clear that $J(A)$ and $J(B)$ are nil. Also, the assumption gives that $A/J(A)$ and $B/J(B)$ are weakly periodic. It now follows from Proposition \ref{proposition1.7} that $A$ and $B$ are weakly periodic, as promised.

($\Leftarrow$). We have $R=S+K$, where
$S=\begin{pmatrix}
     A & 0 \\
     0 & B
\end{pmatrix}$ is a subring of $R$ and
$K=\begin{pmatrix}
   MN & M \\
   N  & NM
\end{pmatrix}$ is a nilpotent ideal of $R$.

Indeed, one computes that $K^{2l}=
\begin{pmatrix}
{(MN)}^l & {(MN)}^lM \\
{(NM)}^lN & {(NM)}^l
\end{pmatrix}$. But one checks that $S\cong A \times B$. Now, the assumption and Proposition \ref{proposition1.6} assure that $S$ is weakly periodic. Finally, Lemma \ref{lemma5.4} ensures the desired result.
\end{proof}

The next three consequences are immediate.

\begin{cor}\label{new}
Let $R$, $S$ be two rings, and let $M$ be an $(R,S)$-bi-module. Then, the formal triangular matrix ring $T(R,S,M)$ is weakly periodic if and only if $R$ and $S$ are weakly periodic rings.
\end{cor}

\begin{cor}\label{corollary5.6}
Let $R$ be a ring and $n\geq 1$ a natural number. Then, $\mathbb{T}_n(R)$ is weakly periodic if and only if $R$ is weakly periodic.
\end{cor}

\begin{cor}\label{corollary5.7}
Let $R$ be a ring with $s\in C(R)\cap \mathrm{Nil}(R)$. Then, the formal matrix ring $K_s(R)$ is weakly periodic if and only if $R$ is weakly periodic.
\end{cor}

\begin{thm}\label{theorem5.8}
Let $R$ be a ring with $s\in C(R)\cap \mathrm{Nil}(R)$ and $n\geq 2$. Then, $\mathbb{M}_n(R;s)$ is weakly periodic if and only if $R$ is weakly periodic.
\end{thm}

\begin{proof}
Arguing as in the proof of Theorem \ref{theorem2.11} and utilizing Theorem \ref{theorem5.5} together with Corollary \ref{corollary5.7}, the claim follows.
\end{proof}

\begin{rem}\label{weakperiod}
In Theorem \ref{theorem5.5}, Corollary \ref{new}, Corollary \ref{corollary5.6}, Corollary \ref{corollary5.7}, and Theorem \ref{theorem5.8} if the requirement "nil-clean" is replaced by the more general "weakly periodic", then the assertions remain true. In this case, Corollary \ref{corollary5.7} and Theorem \ref{theorem5.8} could be considered as extensions of \cite[Corollary 3.5]{KWZ} and \cite[Theorem 3.6]{KWZ}, respectively.
\end{rem}

Some new minor advantage on the Diesl's problem for nil-cleanness of matrix rings over nil-clean rings (see \cite{D13}) is the following one:

\begin{cor}
The formal matrix ring $K_2(R)$ is nil-clean if and only if the ring $R$ is nil-clean.
\end{cor}

\begin{proof}
If $K_2(R)$ is nil-clean, then $2.I_2$ is nilpotent in view of \cite[Proposition 3.14]{D13}. So, one deduces that $2 \in \mathrm{C}(R) \cap \mathrm{Nil}(R)$. If $R$ is nil-clean, then by the cited proposition we again have $2\in \mathrm{C}(R)\cap \mathrm{Nil}(R)$. Hence, in any case, $2 \in \mathrm{C}(R) \cap \mathrm{Nil}(R)$. Now, Corollary~\ref{corollary5.7} and Remark~\ref{weakperiod} apply to get the desired assertion.
\end{proof}

Besides, in the spirit of results from \cite{BCDM} and \cite{KLZ}, one can state the following:

\medskip

\noindent{\bf Conjecture 2:} Suppose $D$ is a division ring and $n\geq 1$ is an integer. Then the matrix ring $\mathbb{M}_n(D)$ is weakly periodic if and only if $D$ is a finite (and hence potent) field.

\medskip

If this statement holds true, then one can extend Theorem~\ref{theorem2.5} to weakly periodic rings, that is, if $R$ is a right (resp., left) perfect weakly periodic ring, then the matrix ring $\mathbb{M}_n(R)$ is weakly periodic. 

Some new relevant comments to these from Remark~\ref{Dis}, which could be some advantage on the currently insurmountable Diesl's problem that we showed above to be somewhat subject to the famous K\"othe's conjecture, are these: It is reasonably logical to ask if $R$ is a (strongly) nil-clean ring and $n\geq 1$ is an arbitrary integer, is then the matrix ring $\mathbb{M}_n(R)$ weakly periodic? In fact, in accordance with \cite{DL} (see also \cite{KWZ}), we may write in the "strongly" case that the factor-ring $R/J(R)$ is boolean and that $J(R)$ is nil. Thus, either \cite{BCDM} or \cite{D2} applies to deduce that $\mathbb{M}_n(R/J(R))\cong \mathbb{M}_n(R)/\mathbb{M}_n(J(R))$ is nil-clean. Now suppose K\"othe's Conjecture is satisfied. Exercise 10.25 in \cite{L} guarantee $J(\mathbb{M}_n(R))=\mathbb{M}_n(J(R))$ is nil. It now follows from \cite[Proposition 3.15]{D13} that $\mathbb{M}_n(R)$ is nil-clean. Therefore, one can say that the solution to the following problem is also subject to K\"othe's Conjecture: {\it If $R$ is a strongly nil-clean ring, then the matrix ring $\mathbb{M}_n(R)$ is nil-clean}.

\medskip

We continue our work with some statements on weakly periodic group rings.

\begin{lem}\label{lemma3.9}
Suppose that the group ring $RG$ is weakly periodic. Then, $R$ is a weakly periodic ring and the center $\mathrm{Z}(G)$ of the group $G$ is a torsion group.
\end{lem}

\begin{proof}
Being a homomorphic image of $RG$, the ring $R$ is too weakly periodic. It follows from Proposition \ref{proposition1.6} (2) that $RG \cong \prod_{i=1}^k R_iG$ such that $p_i \in J(R_i)$, where each $p_i$ is a prime. In virtue of the epimorphism $R_iG \rightarrow (\frac{R_i}{J(R_i)})G$, it is clear that $(\frac{R_i}{J(R_i)})G$ is weakly periodic and $\mathrm{char}(\frac{R_i}{J(R_i)})=p_i$. Thus, without loss of generality, we can assume that the ring $R$ is of prime characteristic $p$.

Suppose now that $g \in \mathrm{Z}(G)$. It is obvious that $g$ is a periodic element. So, $g^n-g \in \mathrm{Nil}(RG)$ for some $n \in \mathbb{N}$, and hence there is an integer $k \geq 2$ such that ${(g^n-g)}^k=0$. That is,  ${(g^{n-1}-1)}^k=0$. Now, we take a natural number $l$ to be big enough such that $p^l \geq k$. That is why, ${(g^{n-1}-1)}^{p^l}=0$. Since $RG$ is of characteristic $p$, we have the equality $g^{(n-1)p^l}=1$. Therefore, $\mathrm{Z}(G)$ is a torsion group, as promised.
\end{proof}

We are now intending to prove our main assertion on weakly periodicity of group rings.

\begin{thm}\label{theorem3.10}
The following two statements are true:
\begin{enumerate}
\item[\rm{(1)}]
Let $R$ be a weakly periodic ring with $p \in \mathrm{Nil}(R)$ and let $G$ be a locally finite $p$-group, where $p$ is a prime. Then, the group ring $RG$ is weakly periodic.
\item[\rm{(2)}]
Let $R$ be a ring and let $G$ be a nilpotent group. If the group ring $RG$ is weakly periodic, then $R$ is a weakly periodic ring and $G$ is a torsion group.
\end{enumerate}
\end{thm}

\begin{proof}
(1) As $G$ is a locally finite $p$-group and $p \in Nil(R)$, it is known by \cite[Proposition 16]{Connell} that the augmentation ideal $\Delta(RG)$ of $RG$ is a nil-ideal. Since $\frac{RG}{\Delta(RG)} \cong R$, it follows from Proposition \ref{proposition1.7} that $RG$ is weakly periodic.

(2) The proof of \cite[Theorem 5.2]{BMA} with a slight change in the argumentation also works in this case.
\end{proof}

We are now ready to prove the aforementioned statement (*) in Section \ref{SectionGroupRing} by utilizing a different approach.

\begin{cor}\label{corollary3.11}
The following two statements are valid:
\begin{enumerate}
\item[\rm{(1)}]
Let $R$ be a periodic ring with $p \in \mathrm{Nil}(R)$ and let $G$ be a locally finite $p$-group, where $p$ is a prime. Then, the group ring $RG$ is periodic.
\item[\rm{(2)}]
Let $R$ be a ring and let $G$ be a group. If the group ring $RG$ is periodic, then $R$ is a periodic ring and $G$ is a torsion group.
\end{enumerate}
\end{cor}

\begin{proof}
(1) It can be established in a way of similarity to item (1) of Theorem \ref{theorem3.10}.

(2) It is evident that $R$ is a periodic ring. Take an arbitrary element $g \in G$. But $R\langle g \rangle$ is a periodic subring of $RG$. Thus, Theorem \ref{theorem3.10} allows us to conclude that $\langle g \rangle$ is a torsion group. Therefore, $G$ is a torsion group too, as stated.
\end{proof}

\section{Open Questions}

In closing, the following six extra problems will, hopefully, motivate a further development of the subject and our possible subsequent research on it.

\medskip

\noindent{\bf Question 1:} If $R$ is a periodic ring, is it true that $\mathbb{M}_n(R)$ is also a periodic ring for any $n\in \mathbb{N}$?

\medskip

If not, we may ask the following two questions (compare with Theorem~\ref{theorem2.3} listed above):

\medskip

\noindent{\bf Question 2:} If $R$ is a periodic ring, does it follow that $\mathbb{M}_n(R)$ is a weakly periodic ring for any $n\in \mathbb{N}$?

\medskip

\noindent{\bf Question 3:} If $R$ is a periodic ring, does it follow that $\mathbb{M}_n(R)$ is a $\pi$-UU ring for any $n\in \mathbb{N}$?

\medskip

If yes, this will contrast \cite[Example 2.9]{CD} where it was established that if $R$ is a $\pi$-UU ring (concretely, the ring of integers $\mathbb{Z}$), then $\mathbb{M}_n(R)$ is {\it not} a periodic ring for any $n\geq 2$.

\medskip

We can also examine the next three queries of interest:

\medskip

\noindent{\bf Question 4:} If $R$ is an UU-ring, is $\mathbb{M}_n(R)$ a $\pi$-UU ring for any $n\in \mathbb{N}$?

\medskip

\noindent{\bf Question 5:} If $R$ is a $\pi$-UU ring, is then $J(R)$ nil?

\medskip

\noindent{\bf Question 6:} If $F$ is a field of prime characteristic $p$ with torsion units and $G$ is a torsion not necessarily Abelian group, is then the group ring $FG$ $\pi$-UU?

\medskip

We just know that if $G$ is an Abelian torsion group, then the group ring $FG$ over such a field $F$ has only torsion units.

\medskip

\noindent{\bf Funding:} The scientific work of the first-named author was supported under the development program of the Volga Region Mathematical Center (agreement no.075-02-2020-1478). The scientific work of the third-named author (Peter V. Danchev) was supported in part by the Bulgarian National Science Fund under Grant KP-06 No. 32/1 of December 07, 2019 as well as by the Junta de Andaluc\'ia, FQM 264.


\medskip


\vskip3.0pc

\begin{thebibliography}{99}


\bibitem{smz-2021}
A.N. Abyzov and D.T. Tapkin, {\it Rings over which every matrices are sums of idempotent and $q$-potent matrices}, Sib. Math. J. \textbf{62} (2021), no. 1, 1--13.

\bibitem{laa-2021}
A.N. Abyzov and D.T. Tapkin, {\it When is every matrix over a ring the sum of two tripotents?}, Lin. Algebra \& Appl. \textbf{630} (2021), 316--325.

\bibitem{AT}
A.N. Abyzov and D.T. Tapkin, {\it On rings with $x^n-x$ nilpotent}, J. Algebra \& Appl. \textbf{21} (2022), no. 6.

\bibitem{A}
S.A. Amitsur, {\it Radicals of polynomial rings}, Canad. J. Math. \textbf{8} (1956), 355--361.

\bibitem{Anderson}
D. D. Anderson and P. V. Danchev, {\it A note on a theorem of Jacobson related to periodic rings}, Proc. Am. Math. Soc. \textbf{148} (2020), no. 12, 5087--5089

\bibitem{Bac95} 
G. Baccella, {\it Semi-artinian $V$-rings and semi-artinian von Neumann regular rings}, J. Algebra \textbf{173} (1995), no. 3, 587--612.

\bibitem{Badawi}
A. Badawi, {\it On Abelian $\pi$-regular rings}, Commun. Algebra 25(4) (1997), 1009--1021. 

\bibitem{BMA}
R. Barati, A. Mousavi and A.N. Abyzov, {\it Rings whose elements are sums of $m$-potents and nilpotents}, Commun. Algebra \textbf{50} (2022), no. 10, 4437--4459.

\bibitem{BM}
R. Barati and A. Mousavi, {\it A note on weakly nil clean rings}, Mediterr. J. Math. \textbf{20} (2023).

\bibitem{Bell}
H.E. Bell and H. Tominaga, {\it On periodic rings and related rings}, Math. J. Okayama Univ. \textbf{28} (1986), 101--103.

\bibitem{BCL}
A.D. Bouzidi, A. Cherchem and A. Leroy, {\it Exponents of skew polynomials over periodic rings}, Commun. Algebra  \textbf{49} (2021), no. 4, 1639--1655.

\bibitem{BCDM}
S. Breaz, G. C\v{a}lug\v{a}reanu, P. Danchev and T. Micu, {\it Nil-clean matrix rings}, Lin. Algebra \& Appl. \textbf{439} (2013), no. 10, 3115--3119.

\bibitem{BDZ}
S. Breaz, P. Danchev and Y. Zhou, {\it Rings in which every element is either a sum or a difference of a nilpotent and an idempotent}, J. Algebra \& Appl. \textbf{15} (2016), no. 8.

\bibitem{wei-xing}
W. Chen and S. Cui, {\it On weakly semicommutative rings}, Commun. Math. Res. \textbf{27} (2011), 179--192.

\bibitem{morphicgrouprings}
J. Chen, Y. Li and Y. Zhou, {\it Morphic group rings}, J. Pure \& Appl. Algebra \textbf{205} (2006), no. 3, 621--639.

\bibitem{CS}
H. Chen and M. Sheibani, Theory of Clean Rings and Matrices, Word Scientific Publishing Company, 2022, 692 pp.
https://doi.org/10.1142/12959

\bibitem{HV}
A.Y.M. Chin and H.V. Chen, {\it On strongly $\pi$-regular group rings}, Southeast Asian Bull. Math. \textbf{26} (2002), no. 3, 387--390.

\bibitem{Connell}
I.G. Connell, {\it On the group rings}, Can. J. Math. \textbf{15} (1963), 650--685.

\bibitem{CD}
J. Cui and P. Danchev, {\it Some new characterizations of periodic rings}, J. Algebra \& Appl. \textbf{19} (2020), no. 12.

\bibitem{PD}
P.V. Danchev, {\it Criteria for unit groups in commutative group rings}, Studia Univ. Babes Bolyai, Math. \textbf{51} (2006), 43--61.

\bibitem{Dnew}
P.V. Danchev, {\it Weakly UU rings}, Tsukuba J. Math. {\bf 40} (2016), no. 1, 101--118.

\bibitem{D1}
P.V. Danchev, {\it On exchange $\pi$-UU unital rings}, Toyama Math. J. \textbf{39} (2017), 1--7.

\bibitem{D2}
P.V. Danchev, {\it Strongly nil-clean corner rings}, Bull. Iran. Math. Soc. \textbf{43} (2017), no. 5, 1333--1339.

\bibitem{DanchevGroupRings}
P.V. Danchev, {\it Commutative periodic group rings}, Mat. Stud. \textbf{53} (2020), 218--220.

\bibitem{D3}
P.V. Danchev, {\it A characterization of weakly tripotent rings}, Rend. Sem. Mat. Univ. Pol. Torino \textbf{79} (2021), no. 1, 21--32.

\bibitem{DGL}
P. Danchev, E. Garc\'ia and M. G\'omez Lozano, {\it Decompositions of matrices into potent and squarezero matrices}, Internat. J. Algebra \& Computat. \textbf{32} (2022), no. 2, 251--263.

\bibitem{DL}
P.V. Danchev and T.-Y. Lam, {\it Rings with unipotent units}, Publ. Math. (Debrecen) \textbf{88} (2016), no. 3-4, 449--466.

\bibitem{DM}
P. Danchev and J. Matczuk, {\it $n$-torsion clean rings}, Contemp. Math. \textbf{727} (2019), 71--82.

\bibitem{DMc}
P.V. Danchev and W.Wm. McGovern, {\it Commutative weakly nil clean unital rings}, J. Algebra {\bf 425} (2015), 410-422.

\bibitem{D13}
A.J. Diesl, {\it Nil clean rings}, J. Algebra \textbf{383} (2013), 197--211.

\bibitem{Di}
A.J. Diesl, {\it Sums of commuting potent and nilpotent elements in rings}, J. Algebra \& Appl. \textbf{22} (2023).

\bibitem{fuchs}
L. Fuchs, Infinite Abelian groups, Volume {\bf I}, Acad. Press, New York and Lodon, 1970.

\bibitem{F}
L. Fuchs, Infinite Abelian Groups, Volume {\bf II}, Acad. Press, New York and London, 1973.

\bibitem{FR}
L. Fuchs and K.M. Rangaswamy, {\it On generalized regular rings}, Math. Z. \textbf{107} (1968), 71--81.

\bibitem{GV}
B. Goldsmith and P. Vamos, {\it A note on clean abelian groups}, Rend. Sem. Mat. Univ. Padova \textbf{117} (2007), 181--191.

\bibitem{G}
K.R. Goodearl, von Neumann Regular Rings, Monographs and Studies in Mathematics, vol. {\bf 4}, Pitman (Advanced Publishing Program), Boston, Massachussetts (1979); 2nd ed., Robert E. Krieger Publishing Co., Inc., Malabar, FL (1991).

\bibitem{wperiodic}
J. Grosen, H. Tominaga and A. Yaqub, {\it On weakly periodic rings, periodic rings and commutativity theorems}, Math. J. Okayama Univ. \textbf{32} (1990), 77--81.

\bibitem{HTY}
Y. Hirano, H. Tominaga and A. Yaqub, {\it On rings in which every element is uniquely expressible as a sum of a nilpotent element and a certain potent element}, Math. J. Okayama Univ. \textbf{30} (1988), 33--40.

\bibitem{Hirano}
Y. Hirano, {\it On periodic P.I. rings and locally finite rings}, Math. J. Okayama Univ. \textbf{33} (1991), 115--120.

\bibitem{HKL}
C. Huh, N.K. Kim and Y. Lee, {\it Examples of strongly $\pi$-regular rings}, J. Pure \& Appl. Algebra \textbf{189} (2004), no. 1-3, 195--210.

\bibitem{NI}
S.U. Hwang, Y.C. Jeon and Y. Lee, {\it Structure and topological conditions of NI rings}, J. Algebra \textbf{302} (2006), 186--199.

\bibitem{kanwar2015clean}
P. Kanwar, A. Leroy and J. Matczuk, {\it Clean elements in polynomial rings}, Contemp. Math. \textbf{634} (2015), 197--204.

\bibitem{khurana2021lifting}
D. Khurana, {\it Lifting potent elements modulo nil ideals}, J. Pure \& Appl. Algebra \textbf{225} (2021), no. 11.

\bibitem{KLZ}
M.T. Ko\c{s}an, T.-K. Lee and Y. Zhou, {\it When is every matrix over a division ring a sum of an idempotent and a nilpotent?}, Lin. Algebra \& Appl. {\bf 450} (2014), no. 11, 7--12.

\bibitem{KWZ}
M.T. Ko\c{s}an, Z. Wang and Y. Zhou, {\it Nil-clean and strongly nil-clean rings}, J. Pure \& Appl. Algebra \textbf{220} (2016), no. 2, 633--646.

\bibitem{KrylovGMR}
P.A. Krylov, {\it Isomorphism of generalized matrix rings}, Algebra \& Logic \textbf{47} (2008), no. 4, 258--262.

\bibitem{krylov2013endomorphism}
P. Krylov, A. Mikhalev and A.A. Tuganbaev, Endomorphism Rings of Abelian Groups, Springer Science \& Business Media (2013).

\bibitem{L}
T.-Y. Lam, A First Course in Noncommutative Rings, Second Edition, Graduate Texts in Math., Vol. \textbf{131}, Springer-Verlag, Berlin-Heidelberg-New York, 2001.

\bibitem{NJ}
C.I. Lee and S.Y. Park, {\it When nilpotents are contained in Jacobson radical}, J. Korean Math. Soc. \textbf{55} (2018), no. 5, 1193--1205.

\bibitem{LeeZhou}
T.-K. Lee and Y. Zhou, {\it Armendariz and reduced rings}, Commun. Algebra \textbf{32} (2004), no. 6, 2287--2299.

\bibitem{Marks}
G. Marks, {\it On $2$-primal Ore extensions}, Commun. Algebra \textbf{29} (2001), no. 5, 2113--2123.

\bibitem{M}
J. Matczuk, {\it Conjugate (nil) clean rings and K\"othe's problem}, J. Algebra \& Appl. \textbf{16} (2017), no. 4.

\bibitem{McGov}
W.Wm. McGovern, {\it The group ring $\mathbb{Z}_{(p)}C_q$ and Ye's theorem}, J. Algebra \& Appl. {\bf 17} (2018), no. 6.

\bibitem{Nasr-Isfahani}
A.R. Nasr-Isfahani, {\it Radicals of skew polynomial and skew Laurent polynomial rings over skew Armendariz rings},
Commun. Algebra \textbf{42} (2014), no. 3, 1337--1349.

\bibitem{N}
W.K. Nicholson, {\it Lifting idempotents and exchange rings}, Trans. Amer. Math. Soc. \textbf{229} (1977), 269--278.

\bibitem{passman}
D.S. Passman, {\it Nil ideals in group rings}, Mich. Math. J. \textbf{9} (1962), no. 4, 375--384.

\bibitem{P}
V. Peri\'c, {\it On rings with polynomial identity $x^n-x=0$}, Publ. Inst. Math. (Beograd) \textbf{34(48)} (1983), 165--167.

\bibitem{pierce1982}
R.S. Pierce, Associative Algebras, Graduate Texts in Math., Springer-Verlag, Berlin-Heidelberg-New York (1982).

\bibitem{PDK}
S. Purkait, T.K. Dutta and S. Kar, {\it On $m$-clean and strongly $m$-clean rings}, Commun. Algebra \textbf{48} (2020), no. 1, 218--227.

\bibitem{S}
A. Stancu, {\it A note on commutative weakly nil clean rings}, J. Algebra \& Appl. \textbf{15} (2016), no. 10.

\bibitem{St}
J. \v{S}ter, {\it Nil-clean quadratic elements}, J. Algebra \& Appl. {\bf 16} (2017), no. 10.

\bibitem{Ster}
J. \v{S}er, {\it On expressing matrices over $\mathbb{Z}_2$ as the sum of an idempotent and a nilpotent}, Linear Algebra \& Appl. \textbf{544} (2018), 339--349.

\bibitem{TZ}
G. Tang and Y. Zhou, {\it A class of formal matrix rings}, Lin. Algebra \& Appl. \textbf{438} (2013), 4672--4688.

\bibitem{TLZ}
G. Tang, C. Li and Y. Zhou, {\it Study of Morita contexts}, Commun. Algebra \textbf{42} (2014), no. 4, 1668--1681.

\bibitem{T}
A. Tuganbaev, Rings Close to Regular, Mathematics and Its Applications, vol. {\bf 545}, Kluwer Academic Publishers, Dordrecht (2002).

\bibitem{W}
T. Wu, {\it On exchange rings with bounded index of nilpotence}, Commun. Algebra {\bf 29} (2001), no. 7, 3089--3098

\bibitem{Ye}
Y. Ye, {\it Semiclean rings}, Commun. Algebra \textbf{31} (2003), no. 11, 5609-5625.

\bibitem{yu1995quasi}
H.P. Yu, {\it On quasi-duo rings}, Glasgow Math. J. \textbf{37} (1995), no. 1, 21--31.

\end{thebibliography}
\end{document}